\documentclass[]{interact}
\usepackage{epstopdf}% To incorporate .eps illustrations using PDFLaTeX, etc.
\usepackage[caption=false]{subfig}% Support for small, `sub' figures and tables

\usepackage[numbers,sort&compress]{natbib}% Citation support using natbib.sty
\bibpunct[, ]{[}{]}{,}{n}{,}{,}% Citation support using natbib.sty
\renewcommand\bibfont{\fontsize{10}{12}\selectfont}% Bibliography support using natbib.sty

\theoremstyle{plain}% Theorem-like structures provided by amsthm.sty

\theoremstyle{definition}

\theoremstyle{remark}

\begin{document}

%\articletype{ARTICLE TEMPLATE}% Specify the article type or omit as appropriate

\title{A Comparative Analysis of Host--Parasitoid Models with Density Dependence Preceding Parasitism}

\author{
\name{Kelsey Marcinko\textsuperscript{a}\thanks{CONTACT Kelsey Marcinko. Email: kmar517@uw.edu} and Mark Kot\textsuperscript{a}}
\affil{\textsuperscript{a}Department of Applied Mathematics, University of Washington, Seattle, WA}
}

\maketitle

\begin{abstract}
We present a systematic comparison and analysis of four discrete-time, host--parasitoid models. For each model, we specify that density-dependent effects occur prior to parasitism in the life cycle of the host. We compare density-dependent growth functions arising from the Beverton--Holt and Ricker maps, as well as parasitism functions assuming either a Poisson or negative binomial distribution for parasitoid attacks. We show that overcompensatory density-dependence leads to period-doubling bifurcations, which may be supercritical or subcritical. Stronger parasitism from the Poisson distribution leads to loss of stability of the coexistence equilibrium through a Neimark--Sacker bifurcation, resulting in population cycles. Our analytic results also revealed dynamics for one of our models that were previously undetected by authors who conducted a numerical investigation. Finally, we emphasize the importance of clearly presenting biological assumptions that are inherent to the structure of a discrete-time model in order to promote communication and broader understanding.
\end{abstract}


%%%%%%%%%%%%%%%%%%%%%%%%%%%%%%%%%%%%%%%%%%%%%%%%%%%%%%%%
\section{Introduction} \label{Section:Intro}
The interactions between insect parasitoids and their hosts are of great interest to ecologists. Roughly 8.5\% of insect species are parasitoids \cite{godfray1994parasitoids}, and they play a significant role in regulating their hosts. Because parasitoid species are specialists on suitable prey, they are often used in biological control programs. This has fueled much interest in developing a better understanding of the dynamics of parasitoids and their hosts. Mathematical models of these host--parasitoid systems are also notable because of the simple and specific modelling assumptions that result from the direct connection between parasitized hosts and parasitoid offspring.

Nicholson and Bailey \cite{nicholson1935balance} laid the foundation for the study of discrete-time host--parasitoid models. Their basic model assumed that oviposition by parasitoids is limited by the number of encounters with hosts and not by parasitoid egg-supply. In addition, they assumed that the number of encounters with hosts is proportional to host abundance and that hosts are equally susceptible to randomly distributed encounters. Their model, however, yields unstable dynamics. As a result, much of the subsequent literature has sought to investigate factors that induce stability. 

In a particularly influential paper, Beddington et al. \cite{beddington1975dynamic} incorporated density-dependent host recruitment, resulting in the model
\begin{subequations}\label{Eq:Beddington}
\begin{align}
N_{t+1}&=N_t e^{r\left(1-\frac{N_t}{K}\right)}e^{-aP_t},\\
P_{t+1}&= c N_t\left(1-e^{-a P_t}\right).
\end{align}
\end{subequations}
Here, $N_t$ is the host density, $P_t$ is the parasitoid density, $r$ is the intrinsic rate of growth, $K$ is the host carrying capacity, $a$ is the parasitoid searching efficiency or area of discovery, and $c$ is the parasitoid clutch size. For a detailed explanation of searching efficiency, see \cite{nicholson1935balance}.

Beddington et al. \cite{beddington1975dynamic} did not specify the life-stage of the host species for which $N_t$ is the density. This is in contrast to Nicholson and Bailey \cite{nicholson1935balance}, who provide extensive biological detail for their model. Beddington et al.'s model also fails to provide a coherent explanation of when the density dependence and parasitism occur during the life-cycle of the host. Specifying the order of events is critical when both density dependence and parasitism affect the host population. 

Model (\ref{Eq:Beddington}) is an example of the more generalized model
\begin{subequations}\label{Eq:MayModel1}
\begin{align}
N_{t+1}&=N_t g(N_t) f(P_t),\\
P_{t+1}&= c N_t[1-f(P_t)].
\end{align}
\end{subequations}
Model (\ref{Eq:MayModel1}) assumes that parasitism affects the original $N_t$ hosts, so that a fraction of hosts, $f(P_t)$, survive parasitism. The survivors then produce offspring with a per capita recruitment, $g(N_t)$, that depends on the original number of hosts. The model also assumes that new parasitoids are produced in proportion to the number of parasitized hosts. Murdoch et al. \cite{murdoch2003consumer} note that the host biology described above is unlikely, though May et. al. \cite{may1981density} provide the example of the winter moth (\emph{Operophtera brumata}) and a fly, \emph{Cyzenis albicans}, that have this biology.

May et al. \cite{may1981density} evaluated model (\ref{Eq:MayModel1}) along with two other models to investigate whether the temporal sequence of host density-dependence and parasitism can affect the dynamics of the populations. In a conclusion that is consistent with the earlier findings of Wang and Gutierrez \cite{wang1980assessment}, May et al. noted that the `sequence of density dependence and parasitism in the host life-cycle can have a significant effect on the population dynamics' \cite{may1981density}. May et al. further recommended that model (\ref{Eq:MayModel1}) be abandoned unless the biology of a particular system demands it.

Numerous investigators  \cite{asheghi2014bifurcations, din2017global, din2016qualitative, hassell1978dynamics, jang2012discrete, kapcak2013stability, kon2006multiple, livadiotis2015discrete} have nevertheless cited Beddington et al. \cite{beddington1975dynamic} and use the structure of model (\ref{Eq:MayModel1}). These authors often derive their models from previous work, without a careful explanation of the underlying biology. Many books  \cite{allen2007introduction, edelstein2005mathematical, murray2003mathematical, turchin2013complex} also present some version of model (\ref{Eq:MayModel1}).  Mills et al \cite{mills1996modelling}, Murdoch et al. \cite{murdoch2003consumer}, and Hassell \cite{hassell2000spatial}, are among the few authors who recognize and discuss the biological assumptions inherent in model (\ref{Eq:MayModel1}).

In this paper, we carefully develop, analyze, and compare four models that assume that density-dependent growth precedes parasitism. These models correspond to a biologically reasonable alternative system presented by May et al. \cite{may1981density}. We consider two functions for density dependence of the host and two functions for parasitism. For each combination of these nonlinear functions, we perform stability analyses to determine dynamics and bifurcations. From these analyses, we conclude that stronger nonlinearity in the density-dependence term produces different effects than stronger parasitism. 

This paper has eight sections. In the second section, we present the biological assumptions underlying our models and the general form of our equations. In the third section, we outline our methods of analysis. In the following four sections, we present four models. The first two models use a fractional function for parasitism, while the next two models use an exponential form. The first and third models assume compensatory density-dependence, while the second and fourth models include overcompensatory density-dependence. The first model yields highly stable dynamics. The second model has a period-doubling route to chaos. A Neimark--Sacker bifurcation occurs in the third model. The fourth model has exponential functions for both density dependence and parasitism, leading to the greatest variability in dynamics. For certain parameter values, there are two interior equilibria, no more than one of which is stable. Both a Neimark--Sacker bifurcation and a subcritical period-doubling bifurcation occur in this model. We conclude with a discussion of the value of understanding the differences between these models. 
%%%%%%%%%%%%%%%%%%%%%%%%%%%%%%%%%%%%%%%%%%%%%%%%%%%%%%%%%%%%%%%%%%%%%%%%%%%%%%%%%%%%%%

\section{Model formulation}\label{Section:Model}
We now consider the model 
\begin{subequations}\label{Eq:Model}
\begin{align}
N_{t+1}&=N_tG(N_t)F(P_t),\\
P_{t+1}&=cN_tG(N_t)[1-F(P_t)].\label{Eq:Model_P}
\end{align}
\end{subequations}
Although this model is consistent with more than one biological scenario, we make several specific choices here. Let $N_t$ be the density of reproducing host adults, and let $P_t$ be the density of adult female parasitoids. $G(N_t)$ is the host per-capita-recruitment. $F(P_t)$, in turn, is the fraction of hosts that escape parasitism, while $1-F(P_t)$ is the fraction of hosts that succumb to parasitism. 

In order to analyze zero-growth isoclines more easily, we let $H(P_t)$ be the fraction of hosts that succumb to parasitism per adult female parasitoid,
\begin{align}\label{eq:Hpt}
H(P_t)= \frac{1-F(P_t)}{P_t}.
\end{align}
System (\ref{Eq:Model}) can now be written
\begin{subequations}\label{Eq:GrowthFirst}
\begin{align}
N_{t+1}&=N_t G(N_t)[1-P_t H(P_t)],\\
P_{t+1}&=c N_t G(N_t)P_t H(P_t),\label{Eq:GrowthFirst_P}
\end{align} 
\end{subequations}
where $c$ is the clutch size. More precisely, $c$ is the average number of female parasitoids laid on a single host that emerge and successfully become reproducing adults. This model is consistent with the second formulation discussed by May et al. \cite{may1981density}.

%%%%%%%%%%%%% FIGURE 1 HERE %%%%%%%%%%%%%%%%
%%%%%%%%%%%%%%%%%%%%%%%%%%%
%%              FIGURE 1                  %%
%%%%%%%%%%%%%%%%%%%%%%%%%%%
\begin{figure}
 \begin{center}
 \includegraphics[width=3in]{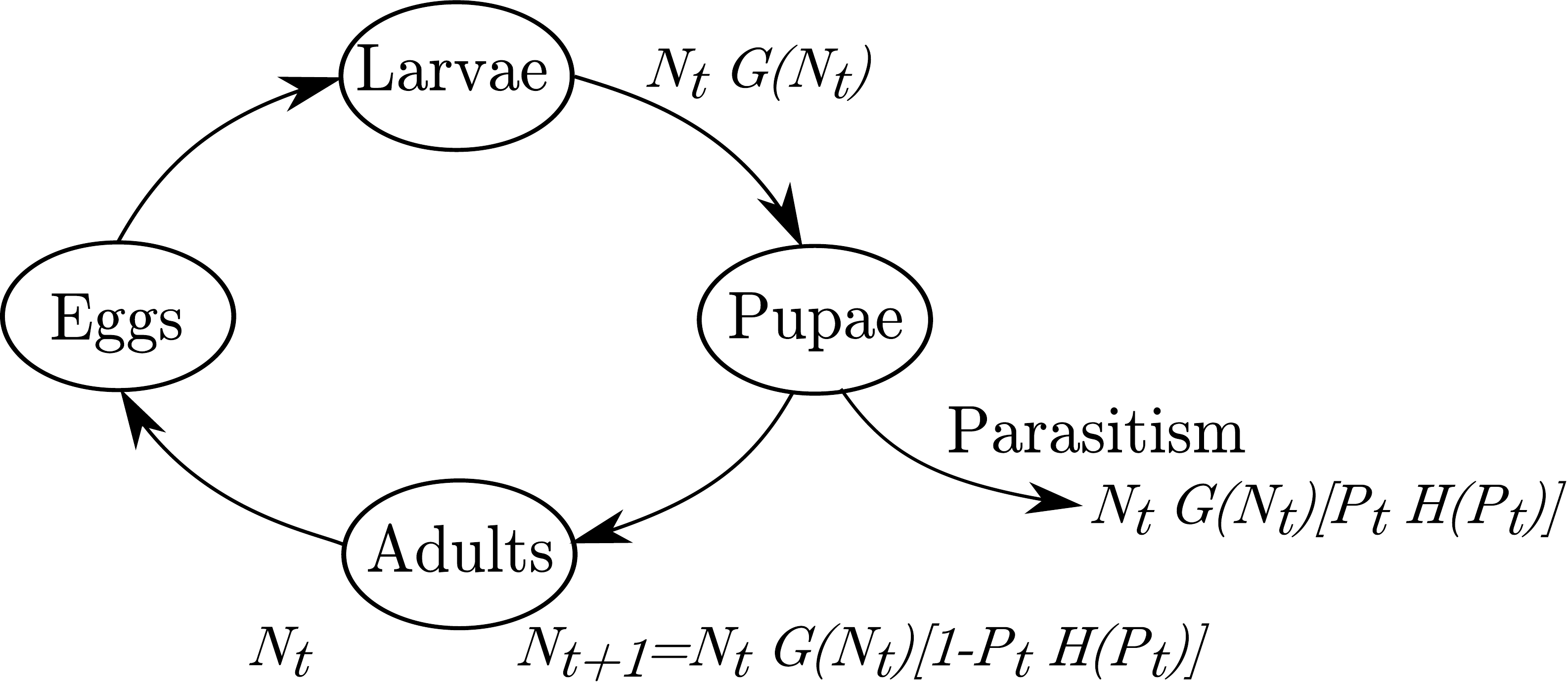} 
\caption{A life-cycle diagram that illustrates a set of biological assumptions that match the formulation of the model set-up with density-dependent competition preceding parasitism. $N_t$ is the density of viable adult hosts that reproduce, and $P_t$ is the density of adult female parasitoids.}\label{fig:LifeCycle}
\end{center}
\end{figure}
%%%%%%%%%%%%%%%%%%

Figure \ref{fig:LifeCycle} illustrates a host life-cycle scenario that matches the biological assumptions of system (\ref{Eq:GrowthFirst}). As above, $N_t$ is the density of reproducing host adults. These adults lay eggs that hatch into larvae. The larvae compete for resources, and $N_t G(N_t)$ larvae survive to the end of larval development. The larvae become pupae, which are parasitized, leaving $N_{t+1}$ adults in the next generation. 

Although the scenario we have described is that of a pupal parasitoid, we emphasize that this is not the only biological scenario described by systems \ref{Eq:Model} and \ref{Eq:GrowthFirst}. The key point, emphasized by Murdoch (\cite{murdoch2003consumer}) and Hassell (\cite{hassell2000spatial}), is that this formulation matches a host life-cycle in which density-dependent competition precedes parasitism. 

We now return to the model. For host density-dependent recruitment, we compare Beverton--Holt growth,
\begin{align}\label{Eq:BevHolt}
N_tG(N_t)&=\frac{R_0N_t}{1+\frac{(R_0-1)}{K}N_t},
\end{align}
and the Ricker curve,
\begin{align}\label{Eq:Ricker}
N_tG(N_t)&=N_te^{r\left(1-\frac{N_t}{K}\right)}.%=N_tR_0^{\left(1-\frac{N_t}{K}\right)}.
\end{align}
Here $R_0=\exp(r)$ is the net reproductive rate, $r=\ln(R_0)$ is the intrinsic rate of growth, and $K$ is the carrying capacity. Recall that the Beverton--Holt growth function is compensatory while the Ricker growth function is overcompensatory.

Early investigators \cite{thompson1924, nicholson1935balance} used the zero term of the Poisson distribution for $F(P_t)$, the fraction of hosts that escape parasitism. May \cite{may1978host} considered varying levels of aggregation and proposed the use of the zero term of the negative binomial distribution,
\begin{align}\label{Eq:NegBinom}
F(P_t)&=\left(1+\frac{aP_t}{\kappa}\right)^{-\kappa}.
\end{align} 
May's use of this function influenced Livadiotis et al. \cite{livadiotis2016kappa}, who studied system (\ref{Eq:MayModel1}) with $\kappa$-parameterized functions for both parasitism and density-dependent intraspecific competition. The formulation used by Livadiotis et al. highlights the similarities in the exponential ($\kappa \to\infty$) and rational ($\kappa = 1$) functions most commonly used for $F(P_t)$ and $G(N_t)$.

In this paper, we focus on two forms of May's function, given by $\kappa=1$ and $\kappa\to \infty$. From equation (\ref{eq:Hpt}), these values of $\kappa$ give the fraction of hosts that succumb to parasitism per adult female parasitoid as
\begin{align}\label{Eq:FracPar}
H(P_t)&=\frac{1}{1+aP_t}
\end{align}
and
\begin{align}\label{Eq:ExpPar}
H(P_t)&=\frac{1}{P_t}\left(1-e^{-aP_t}\right)
\end{align} 
respectively.

Other than in May et al.'s paper \cite{may1981density}, system (\ref{Eq:GrowthFirst}) has not been studied in a way that compares functional forms for modelling parasitism and density dependence. Using equations (\ref{Eq:BevHolt}), (\ref{Eq:Ricker}), (\ref{Eq:FracPar}), and (\ref{Eq:ExpPar}), we will formulate four possible models and compare their dynamics in Sections \ref{Section:M1}--\ref{Section:M4}.

%%%%%%%%%%%%%%%%%%%%%%%%%%%%%%%%%%%%%%%%%%%%%%%%%%%%%%%%%%%%%%%%%%%%%%%%%%%%%%%%%%%%%%%%
\section{Methods of analysis}\label{sec:MoA}
Each of our four models can be written in the general ``density-dependence first" form of system (\ref{Eq:GrowthFirst}). We now nondimensionalize. If we let $y_t=aP_t$, $x_t=N_t/K$, and $b=acK$, we obtain
\begin{subequations}\label{Eq:uvSystem}
\begin{align}
x_{t+1}&=x_t u(x_t,y_t),\\
y_{t+1}&=y_t v(x_t,y_t),
\end{align}
\end{subequations}
where
\begin{align}
u(x_t,y_t)&=g(x_t)[1-y_th(y_t)], \label{Eq:u}\\
v(x_t,y_t)&=bx_tg(x_t)h(y_t)\label{Eq:v}.
\end{align}

For Beverton-Holt growth,
\begin{align}\label{Eq:BevHoltxy}
g(x_t)&=\frac{R_0}{1+(R_0-1)x_t},
\end{align}
 while for Ricker growth, 
\begin{align}\label{Eq:Rickerxy}
g(x_t)&=e^{r(1-x_t)}.
\end{align}
We will call (\ref{Eq:BevHoltxy}) \emph{fractional per-capita-recruitment}, which produces compensatory density-dependence, and (\ref{Eq:Rickerxy})  \emph{exponential per-capita-recruitment}, which produces overcompensatory density-dependence. 

Similarly, the fraction of hosts that succumb to parasitism, $yh(y)$, can be rewritten with
\begin{align}
h(y_t)&=\frac{1}{1+y_t},\label{Eq:fracxy}
\end{align} 
for $\kappa=1$, and
\begin{align}\label{Eq:expxy}
h(y_t)&=\frac{1}{y_t}\left(1-e^{-y_t}\right),
\end{align}
for $\kappa\to\infty$. We will refer to (\ref{Eq:fracxy}) as \emph{fractional parasitism} and (\ref{Eq:expxy}) as \emph{exponential parasitism}.

\begin{figure}
\centering
\subfloat[Host per capita recruitment]{%
       \includegraphics[width = 0.3\linewidth]{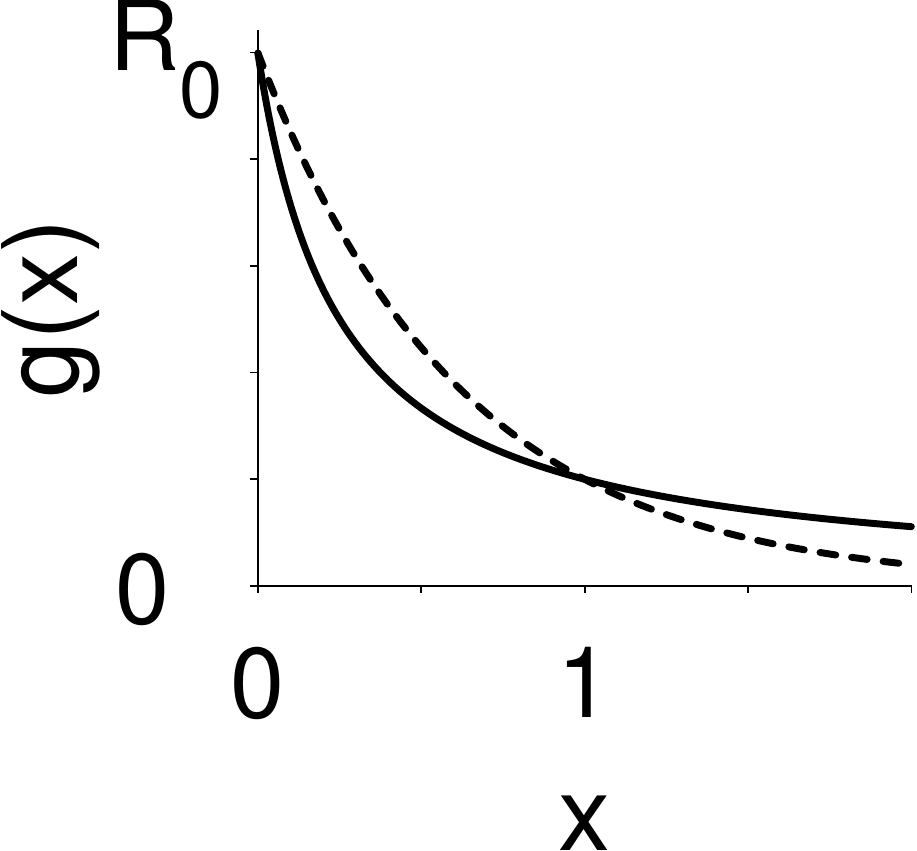}\label{fig:g_x_plot}
     }
\hfill
\subfloat[Host recruitment]{%
       \includegraphics[width = 0.3\linewidth]{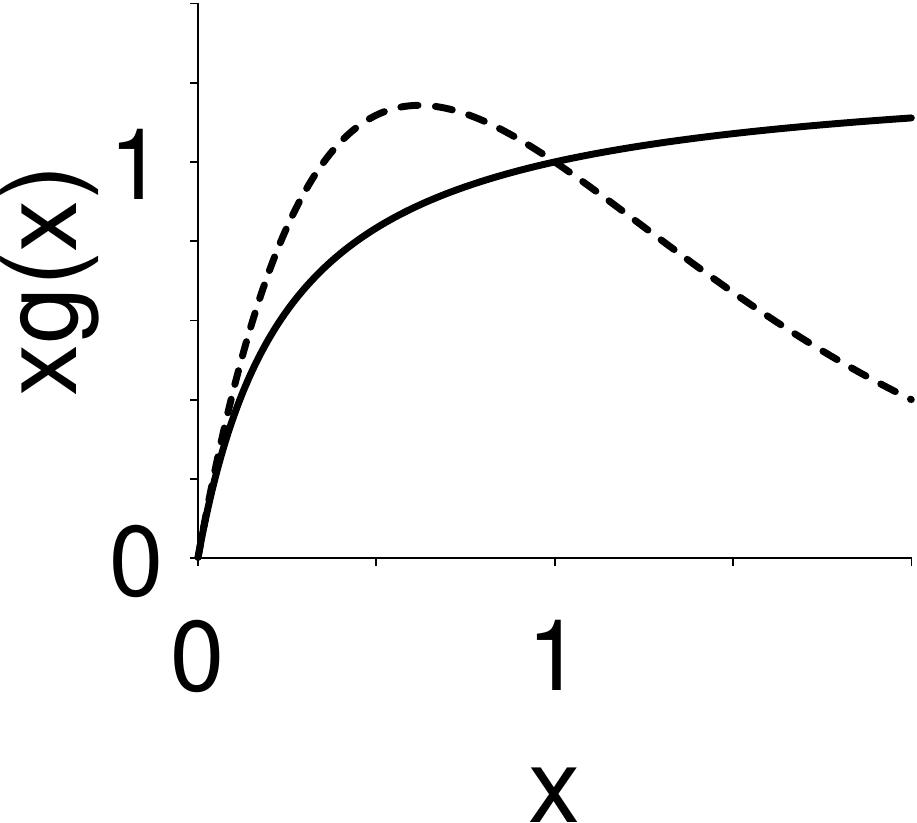}\label{fig:x_g_x_plot}
     }
\hfill
\subfloat[Parasitism]{%
       \includegraphics[width =  0.3\linewidth]{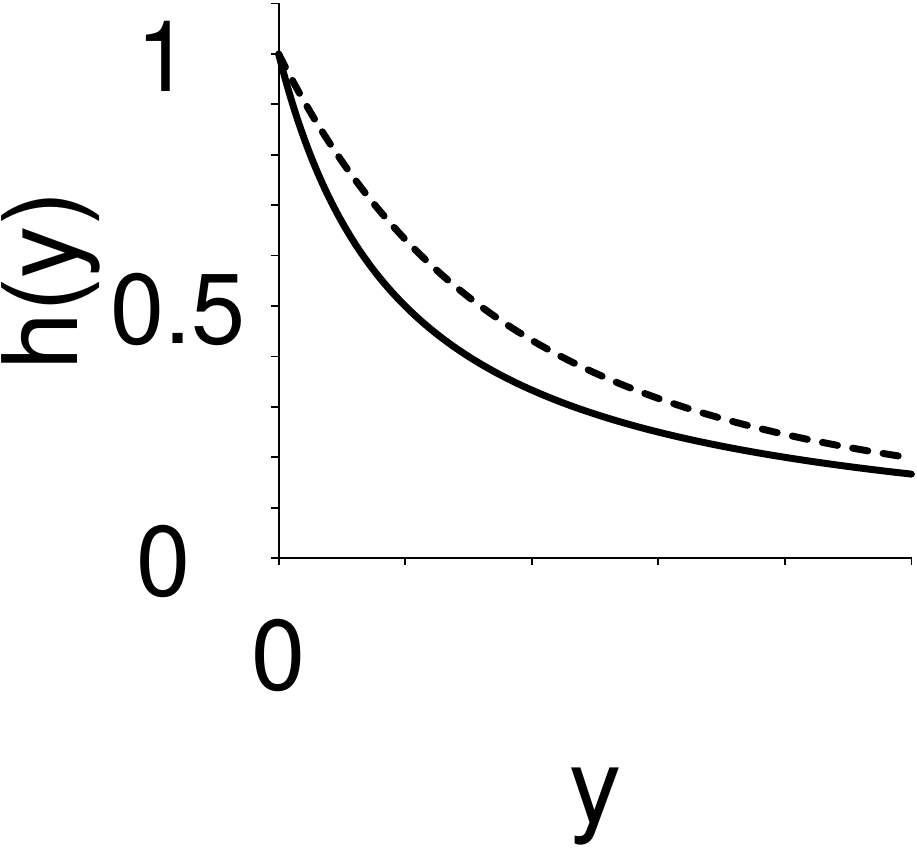}\label{fig:h_y_plot}
     }
 \caption{These figures illustrate the behavior of $g(x)$, $xg(x)$, and $h(y)$ for functions used in our models. Fractional forms of $g(x)$ and $h(y)$ from equations (\ref{Eq:BevHoltxy}) and (\ref{Eq:fracxy}) are shown with solid lines. Exponential forms of $g(x)$ and $h(y)$ from equations (\ref{Eq:Rickerxy}) and (\ref{Eq:expxy}) are shown with dashed lines. Both functions for $g(x)$ are monotonically decreasing from $R_0$. Recruitment, $xg(x)$, is non-monotonic for the exponential form, while it is monotonic for the fractional form. Both fractional and exponential forms of $h(y)$ are positive and monotonically decreasing. The dashed curve remains above the solid curve as $y$ increases.}\label{fig:functionplots}
\end{figure}

In all that follows, we assume $R_0\geq 1$ ($r\geq 0$), since we choose to consider cases where the host species can persist in the absence of the parasitoid species. For $R_0>1$ ($r>0$), the per-capita recruitment, $g(x_t)$, is a positive, monotonically decreasing function that starts from $R_0=\ln(r)$ at $x_t=0$ and crosses 1 at $x_t=1$. Similarly, $h(y_t)$ is positive and monotonically decreasing, with $h(0)=1$. Sample plots of $g(x)$, $xg(x)$, and $h(y)$ are shown in Figure \ref{fig:functionplots}.

%%%%%%%%%% FIGURE 2 HERE %%%%%%%%%%%%%%%%%%%%%%%%%%%
%%%%%%%%%%%%%%%%%%%%%
%%         FIGURE 2          %%%%%%%%%%
%%%%%%%%%%%%%%%%%%%%%%
%%%%%%%%%%%%%%%%%%%%%%%%%%%%%%%%%%%%%%%%%%%%%%%

To find the equilibria of system (\ref{Eq:uvSystem}), we set $x_{t+1}=x_t$ and $y_{t+1}=y_t$. The equilibria occur at (0,0), (1,0), and at solutions of the system
\begin{subequations}\label{Sys:uv1}
\begin{align}
1&=u(x,y)=g(x)[1-yh(y)], \label{Eq:uis1}\\
1&=v(x,y)= bxg(x)h(y),\label{Eq:vis1}
\end{align}
\end{subequations}
where we drop the $t$ subscripts for notational simplicity. For each of our models, it can be shown that $b>1$ is a necessary and sufficient condition for the existence of a \emph{unique} positive solution to system (\ref{Sys:uv1}). For the fourth model, there is a region below $b=1$ for which two positive solutions to system (\ref{Sys:uv1}) exist.

To determine the stability of the equilibria, we form the Jacobian matrix of partial derivatives for system (\ref{Eq:uvSystem}),
\begin{equation}\label{Eq:Jacobian}
J(x,y)=\begin{pmatrix}
xu_x+u & xu_y\\
yv_x & yv_y+v
\end{pmatrix}.
\end{equation}
After evaluating the partial derivatives, the Jacobian may be rewritten
\begin{equation}
J(x,y)=\begin{pmatrix}
[x g'(x) + g(x)][1-y h(y)] & -x g(x)[y h'(y)+h(y)] \\
 & \\
b y h(y)[x g'(x) + g(x)] & b x g(x)[y h'(y) + h(y)]
\end{pmatrix},
\end{equation} 
where we factor to separate the $x$ and $y$ dependencies. We now use the Jacobian evaluated at each of the equilibria to determine stability. 

\subsection{Extinction equilibrium}
At the \emph{extinction point} $(0,0)$, the Jacobian,
\begin{align}
J(0,0) &= \begin{pmatrix}
g(0) & 0 \\
0 & 0
\end{pmatrix}
= \begin{pmatrix}
R_0 & 0 \\
0 & 0
\end{pmatrix},
\end{align}
has eigenvalues $R_0$ and 0. Note that we used $g(0)=R_0$, which was mentioned previously. The extinction equilibrium is unstable for $R_0> 1$. The zero eigenvalue indicates that for initial conditions with $x=0,\ y>0$, the system will collapse to the $(0,0)$ fixed point at the next generation due to the lack of hosts.

\subsection{Exclusion equilibrium}
The equilibrium point (1,0) is known as an \emph{exclusion point} \cite{asheghi2014bifurcations, kapcak2013stability, kapcak2016stability}. Here, the host population persists at carrying capacity, while the parasitoid population goes extinct. The Jacobian for this system is
\begin{align}
J(1,0) &= \begin{pmatrix}
g'(1) +g(1)\ \ \  & -h(0) \\
 & \\
0 & b h(0)
\end{pmatrix}
  = \begin{pmatrix}
 g'(1) + 1 \ \ \ & -1 \\
 & \\
0 & b
\end{pmatrix},
\end{align}
since $h(0) = 1$ for equations (\ref{Eq:fracxy}) and (\ref{Eq:expxy}). The eigenvalues for this triangular system are thus
\begin{equation} \label{Eq:exclEigs}
\lambda_1 = g'(1) + 1, \quad \lambda_2 = b.
\end{equation}

Recall that $g'(1)$ is negative since $g(x)$ is monotone decreasing for $R_0>1$ ($r>0$). Based on the eigenvalues in (\ref{Eq:exclEigs}), we conclude that we need both
\begin{align} \label{Eq:lambda1}
-2<g'(1)<0
\end{align}
and $b<1$ for the exclusion equilibrium to be asymptotically stable. The second inequality in condition (\ref{Eq:lambda1}) is satisfied, so we will check the first inequality for both forms of the host per-capita-recruitment, $g(x)$. 

For equation (\ref{Eq:BevHoltxy}),
\begin{align}
g'(1)&=\frac{(1-R_0)}{R_0},
\end{align}
and the first inequality in (\ref{Eq:lambda1}) becomes
\begin{equation}
-2R_0 < 1-R_0,
\end{equation}
which simplifies to $-1<R_0$. Since the net reproductive rate, $R_0$, is positive, this inequality is true, and the stability of the exclusion equilibrium point hinges on the value of $b$ for our models that use fractional recruitment. For $b<1$, the equilibrium is asymptotically stable, and for $b>1$, the equilibrium is unstable.

For equation (\ref{Eq:Rickerxy}),
\begin{align}
g'(1)&=-r.
\end{align}
Stability thus requires $-2<-r<0$. So for our models that use exponential recruitment, both $b<1$ and $0<r<2$ are necessary for asymptotic stability of the exclusion equilibrium.

\subsection{Coexistence equilibria} \label{Sec:CoexEq}
 The \emph{coexistence equilibria} are the solutions to system (\ref{Sys:uv1}). Biologically, coexistence occurs when both $x$ and $y$ are positive. These equilibria can be explicitly determined for models using fractional parasitism, but not for exponential parasitism. Nevertheless, the coexistence equilibria can be approximated numerically for all cases.

Using equations (\ref{Eq:uis1}) and (\ref{Eq:vis1}), Jacobian matrix (\ref{Eq:Jacobian}) simplifies to
\begin{align}
J(x,y)&=\begin{pmatrix}\label{Eq:JacobianSimp}
xu_x+1 & xu_y\\
yv_x & yv_y+1
\end{pmatrix}.
\end{align} To avoid unnecessarily complicated algebra, we will not proceed from eigenvalues. 

Instead, to determine the conditions for asymptotic stability of the coexistence equilibria, we will apply the Jury conditions \cite{jury1964theory} to each model. These necessary and sufficient conditions for asymptotic stability are
\begin{gather}
1-\tau+\Delta>0, \label{Eq:Jury1} \\
1+\tau+\Delta>0, \label{Eq:Jury2}\\
\Delta<1, \label{Eq:Jury3}
\end{gather}
where $\tau$ is the trace and $\Delta$ is the determinant of the Jacobian matrix evaluated at the implicit or explicit coexistence equilibrium. For matrix (\ref{Eq:JacobianSimp}),
\begin{align}\label{Eq:tau}
\tau &= 2 +xu_x+yv_y,
\end{align}
and
\begin{equation}\label{DetExp}
\begin{split}
\Delta &=1+xu_x+yv_y+xy(u_xv_y-u_yv_x).
\end{split}
\end{equation}

Using these expressions, the first Jury condition, inequality (\ref{Eq:Jury1}), simplifies to
\begin{equation}\label{Jury1}
xy(u_xv_y-u_yv_x)>0.
\end{equation} The first Jury condition will be violated for parameter values such that $x=0$ or $y=0$. For a true coexistence equilibrium point with positive $x$ and $y$ values, inequality (\ref{Jury1}) requires
\begin{align} \label{Jury1part}
u_xv_y-u_yv_x&>0.
\end{align}

We now consider the second Jury condition, inequality (\ref{Eq:Jury2}).  After we write the inequality in terms of $u,v,x,$ and $y$, the condition simplifies to
\begin{equation}
4+2xu_x+2yv_y+xy(u_xv_y-u_yv_x)>0. \label{Jury2}
\end{equation}

Finally, the third Jury condition, inequality (\ref{Eq:Jury3}), can be expressed as
\begin{equation}
1+xu_x+yv_y+xy(u_xv_y-u_yv_x)<1. \label{Jury3}
\end{equation}
These three Jury conditions (\ref{Jury1}, \ref{Jury2}-\ref{Jury3}) will be used for each specific model to determine the requirements on parameters $b$ and $R_0$ (or $r$) to ensure that the coexistence equilibrium is stable.
%%%%%%%%%%%%%%%%%%%%%%%%%%%%%%%%%%%%%%%%%%%%%%%%%%%%
\section{Model 1: Compensatory host density-dependence and fractional parasitism}\label{Section:M1}
The first model we consider uses fractional per-capita-recruitment (\ref{Eq:BevHoltxy}) and fractional parasitism (\ref{Eq:fracxy}). The model is thus
\begin{subequations}\label{Eq:Model1}
\begin{align}
x_{t+1}&=\left[\frac{R_0x_t}{1+(R_0-1)x_t}\right]\left(\frac{1}{1+y_t}\right),\\
y_{t+1}&=b\left[\frac{R_0x_t}{1+(R_0-1)x_t}\right]\left(\frac{y_t}{1+y_t}\right).
\end{align}
\end{subequations}

The coexistence equilibrium for this system is
\begin{equation}\label{Eq:M1Equil}
(x^*,y^*)=\left(\frac{1}{b},\frac{R_0}{1+\left(R_0-1\right)\left(\frac{1}{b}\right)}-1\right).
\end{equation}
 As shown in Appendix \ref{app:Model1FirstQuadrant}, for $R_0>1$, this equilibrium is in the interior of the first quadrant if and only if $b>1$. For $b=1$, the equilibrium given by equation (\ref{Eq:M1Equil}) is the exclusion equilibrium, $(1,0)$.  For $R_0=1$, system (\ref{Eq:Model1}) has a line of equilibria on the $x$-axis, and (\ref{Eq:M1Equil}) reduces to $(1/b,0)$.

 \subsection{Stability region}
\begin{figure}
\centering
\subfloat[Model 1]{%
       \includegraphics[width = 0.45\linewidth]{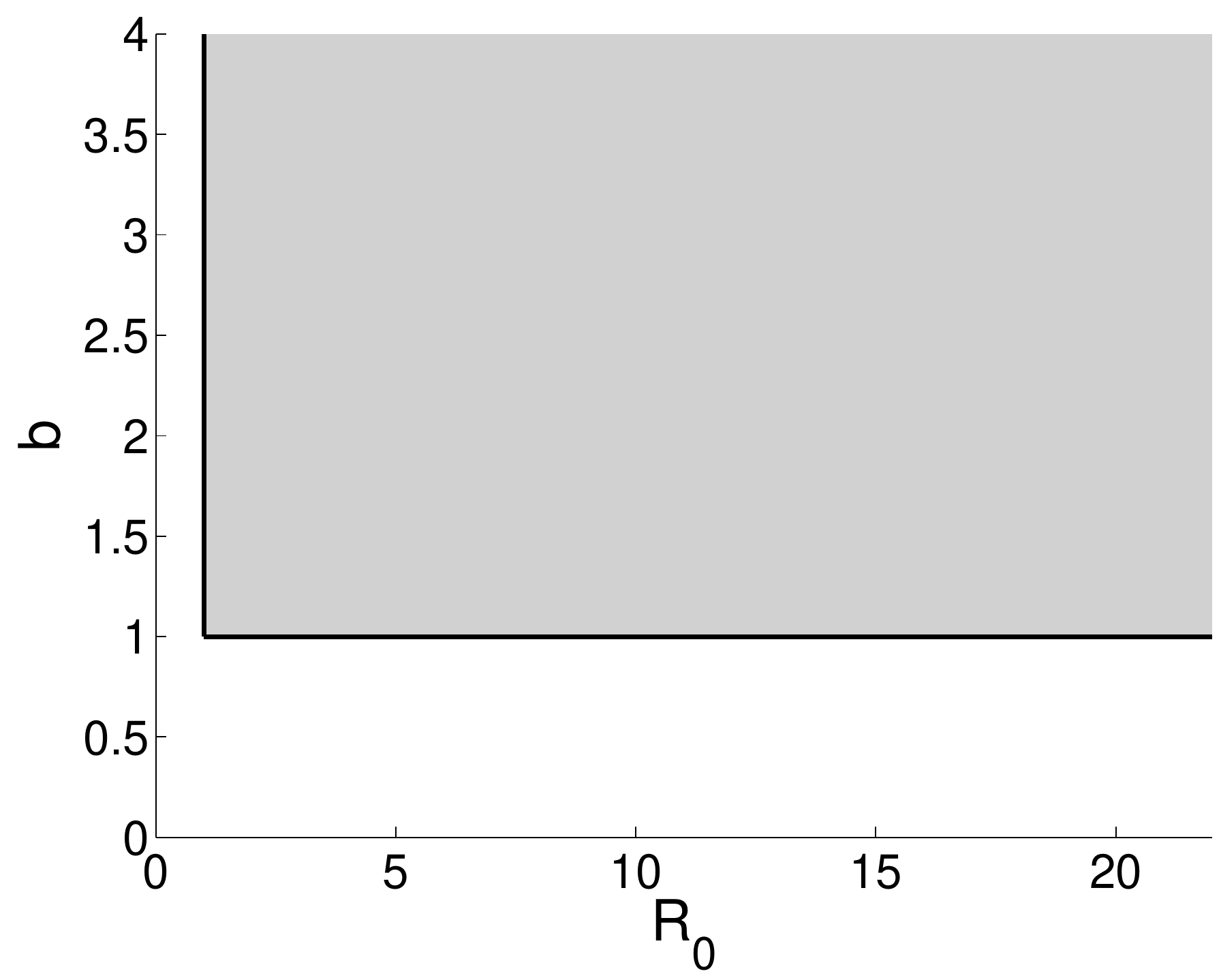}\label{fig:Model1StabilityRegion}
     }
\hfill
\subfloat[Model 2]{%
       \includegraphics[width =  0.45\linewidth]{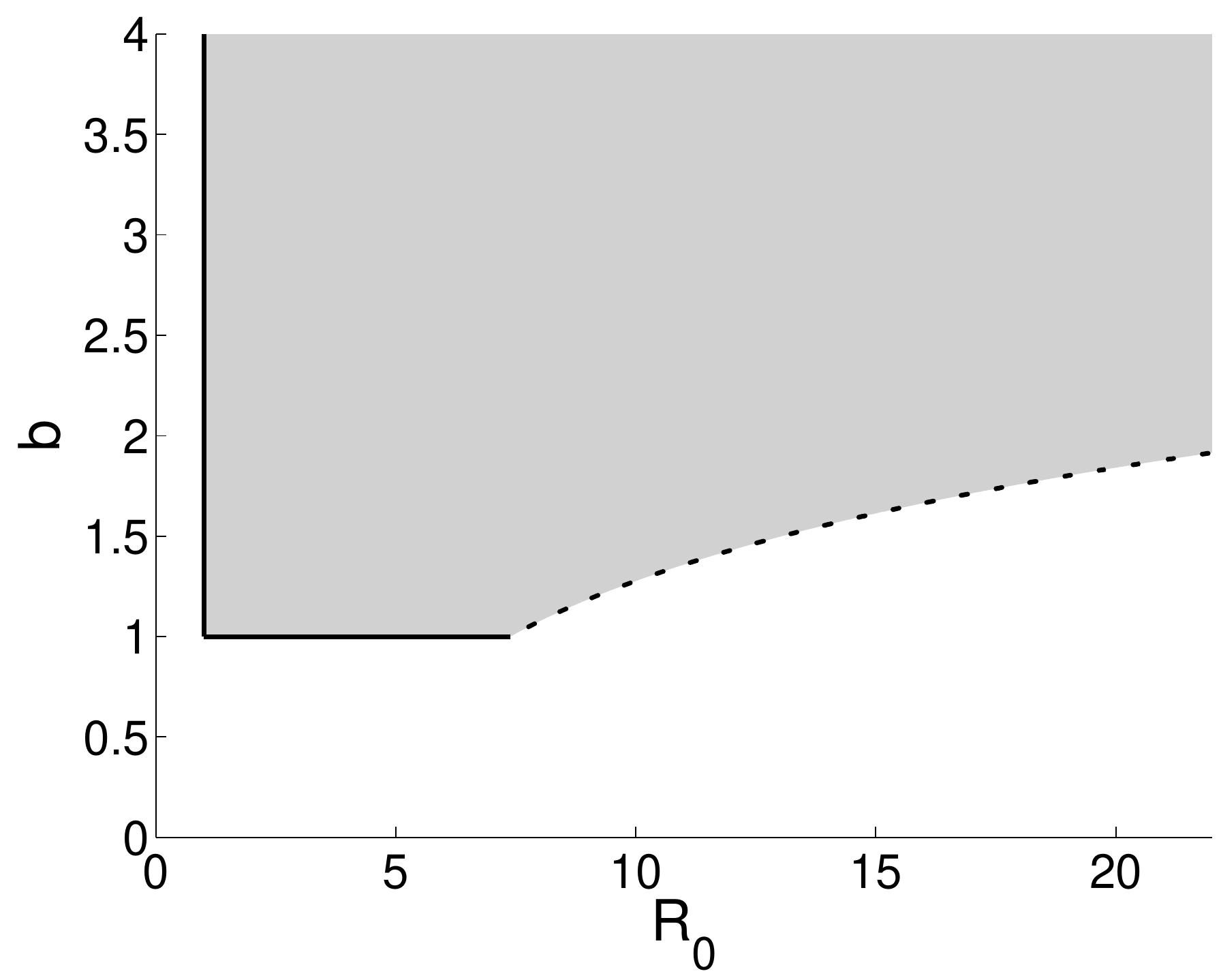}\label{fig:Model2StabilityRegion}
     }
\\
     \subfloat[Model 3]{%
       \includegraphics[width = 0.45\linewidth]{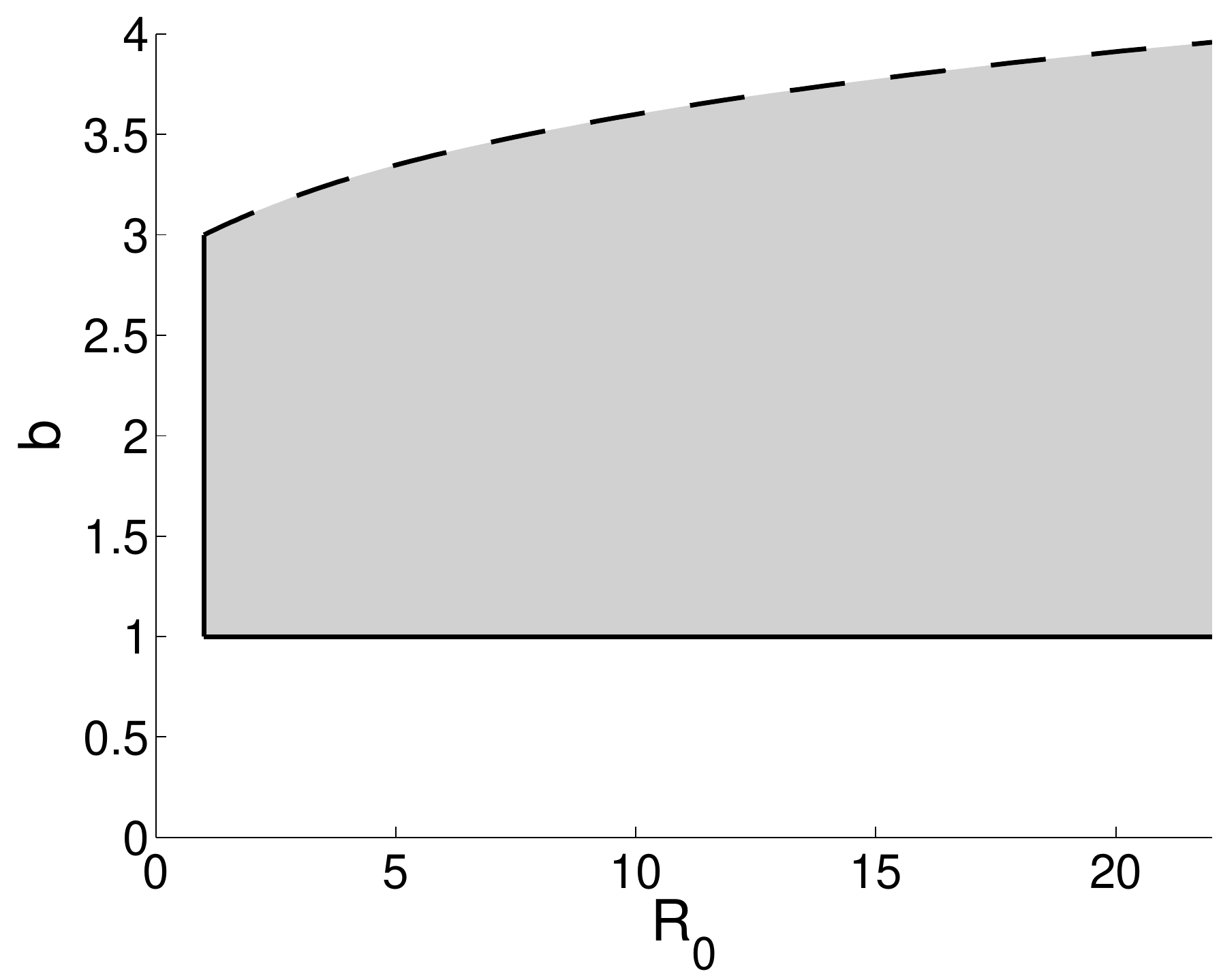}\label{fig:Model3StabilityRegion}
     }
\hfill
\subfloat[Model 4]{%
       \includegraphics[width = 0.45\linewidth]{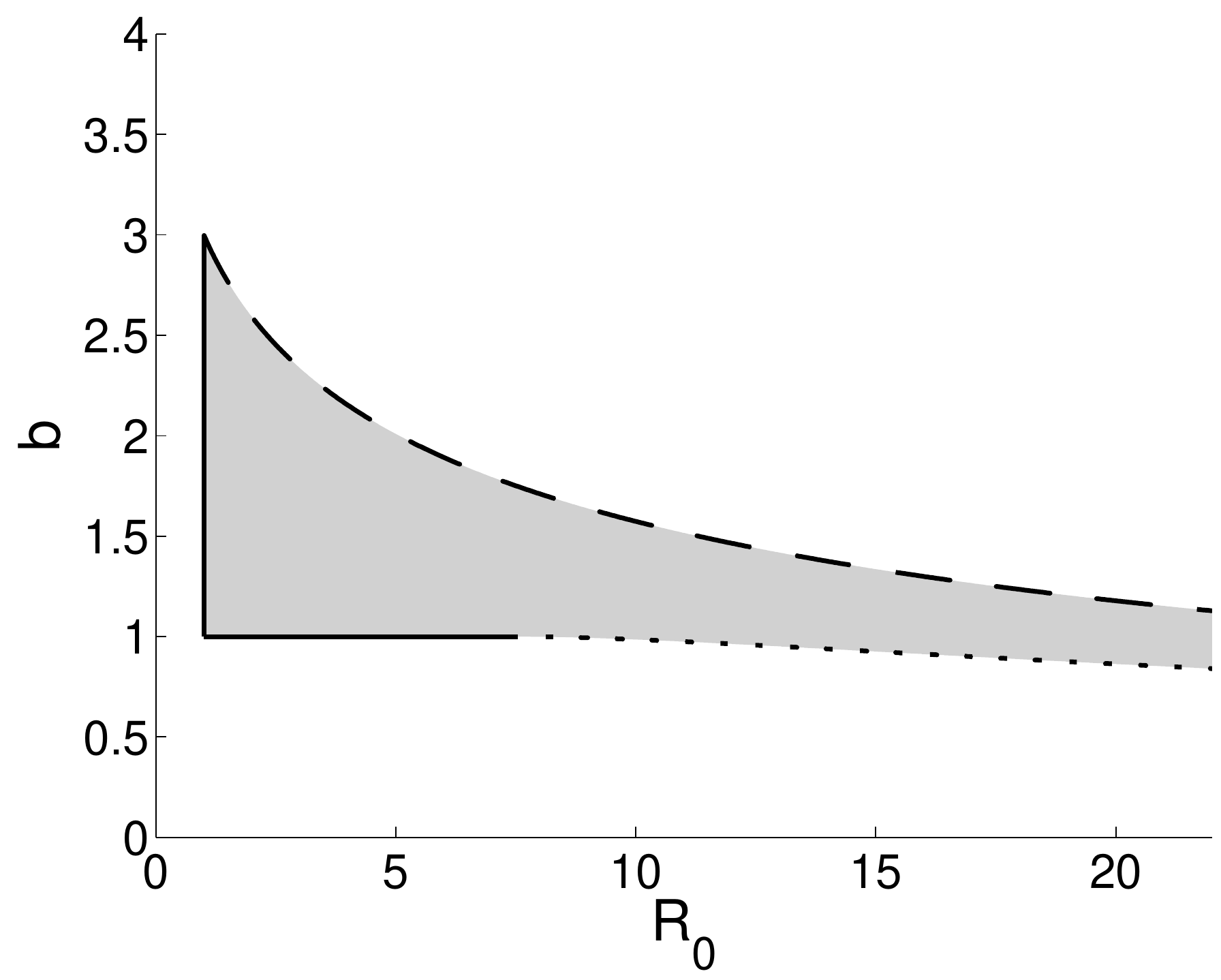}\label{fig:Model4StabilityRegion}
     }
 \caption{Stability regions for the positive coexistence equilibrium for Models 1--4. For $R_0 > 1$, the condition $b>1$ ensures that the coexistence equilibrium is in the interior of the first quadrant for Models 1--3. For Model 4, the system has a unique positive coexistence equilibrium for $b>1$. For $R_0>e^2$, there is a region below $b=1$ for which there are two positive coexistence equilibria. The coexistence equilibrium with the larger $y$ coordinate is stable in the shaded region below $b=1$. For all models, $R_0>1$ is necessary for asymptotic stability. The vertical solid lines are boundary curves where \emph{both} the first and third Jury conditions are violated. The horizontal solid lines are boundary curves where the first Jury condition is violated. The dotted curves are the boundaries where the second Jury condition is violated. The dashed curves are the boundaries where only the third Jury condition is violated. }
\end{figure}

Compensatory (fractional) host recruitment, $xg(x)$, and fractional parasitism are both rational functions, which correspond to a low $\kappa$ index in the parameterized families of common recruitment and parasitism functions (see Livadiotis et al. \cite{livadiotis2016kappa}). When we use fractional per-capita-recruitment for $g(x)$ and fractional parasitism for $h(y)$, the model has a large stability region as seen in Figure \ref{fig:Model1StabilityRegion}. All three Jury conditions are satisfied for the region in parameter space defined by $b>1$, $R_0>1$.  The first Jury condition is violated for $b=1$. Crossing this line corresponds to a transcritical bifurcation. Both the first and third Jury conditions are violated for $R_0=1$. Details are given in Appendix \ref{app:Model1Jury1}--\ref{app:Model1Jury3}. Satisfying the three Jury conditions ensures that the coexistence equilibrium is asymptotically stable.

%%%%%%%% FIGURE 3 HERE %%%%%%%%%%%%%%%%%%
%%%%%%%%%%%%%%%%%%%%%%%%%%%%%
%%       FIGURE 3                   %%%%%%%%%%
%%%%%%%%%%%%%%%%%%%%%%%%
%%%%%%%%%%%%%%%%%%%%%%%%%%%%%%%%%%%%%%%%%%%%%%%
%%%%%%%%%%%%%%%%%%%%%%%%%%%%%%%%%%%%%%%%%%%%%%%%%%

%%%%%%%%%%%%%%%%%%%%%%%%%%%%%%%%%%%%%%%%%%%
\section{Model 2: Overcompensatory host density-dependence and fractional parasitism}\label{Section:M2}
Our second model also uses fractional parasitism, but it incorporates the exponential per-capita-recruitment from equation (\ref{Eq:Rickerxy}). As seen in Figure (\ref{fig:x_g_x_plot}), exponential recruitment is nonmonotonic, so we have introduced stronger nonlinearity in the density-dependence term. These choices yield the model
\begin{subequations}\label{Eq:Model2}
\begin{align}
x_{t+1}&=x_t e^{r(1-x_t)}\left(\frac{1}{1+y_t}\right),\\
y_{t+1} &= b x_t e^{r(1-x_t)}\left(\frac{y_t}{1+y_t}\right).
\end{align}
\end{subequations}

The coexistence equilibrium for this system,
\begin{equation}\label{Eq:M2Equil}
(x^*,y^*)=\left(\frac{1}{b}, e^{r\left(1-1/b\right)}-1\right),
\end{equation}
 is again in the interior of the first quadrant if $r>0$, $b>1$. This is shown in Section \ref{app:Model2FirstQuadrant}. For $b=1$, the equilibrium given by equation (\ref{Eq:M2Equil}) is the exclusion equilibrium, $(1,0)$.  For $r=0$, system (\ref{Eq:Model2}) has a line of equilibria on the $x$-axis, and (\ref{Eq:M2Equil}) reduces to $(1/b,0)$. 
 
\subsection{Stability region}
 As was true for Model 1, the first Jury condition is satisfied for $b>1$, $r>0$ ($R_0>1$), and the third Jury condition is satisfied for $r>0$ ($R_0>1$). Substituting the exponential form of density dependence in place of the fractional form from Model 1 introduces an additional stability criterion for the coexistence equilibrium for Model 2. The second Jury condition is now satisfied above the curve defined, for $u > 3/2$, by
\begin{equation}\label{Eq:M2J2}
r = u-\ln\left(2u-3\right),\quad b=1-\frac{1}{u}\ln(2u-3).
\end{equation}
The derivation of these stability criteria is shown in Sections \ref{app:Model2Jury1}--\ref{app:Model2Jury3}.

The stability region for the coexistence equilibrium is shown in Figure \ref{fig:Model2StabilityRegion}. Crossing the line, $b=1$, $1<R_0<e^2$ violates the first Jury condition, resulting in a transcritical bifurcation. Crossing the line $R_0=1$ violates both the first and third Jury conditions. Crossing the dotted curve from the left in Figure \ref{fig:Model2StabilityRegion} means that one of the real eigenvalues exceeds $-1$ in magnitude. This corresponds to a period-doubling or flip bifurcation. However, the stability analysis holds only in a neighborhood of the equilibrium point. We discuss below the existence of other stable phenomena for this model, including 2-cycles, 4-cycles, and invariant circles. 

\subsection{Bifurcations and attractors}
\begin{figure}
\centering
\subfloat[$b=1.1$]{%
       \includegraphics[width =  0.48\linewidth]{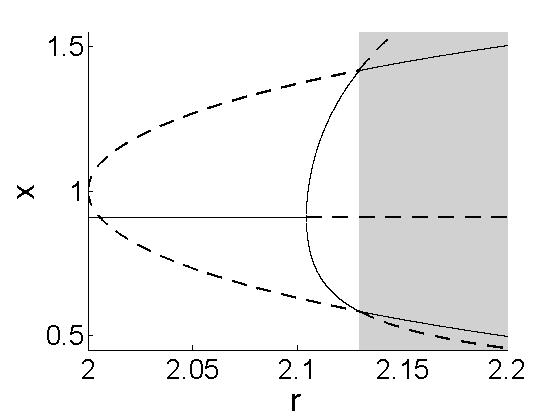}
     }
\hfill
\subfloat[$b=1.1$]{%
       \includegraphics[width = 0.48\linewidth]{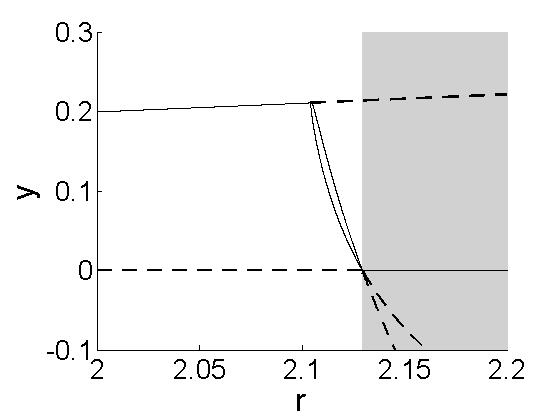}
     }
\\
\subfloat[$b=1.3$]{%
       \includegraphics[width =  0.48\linewidth]{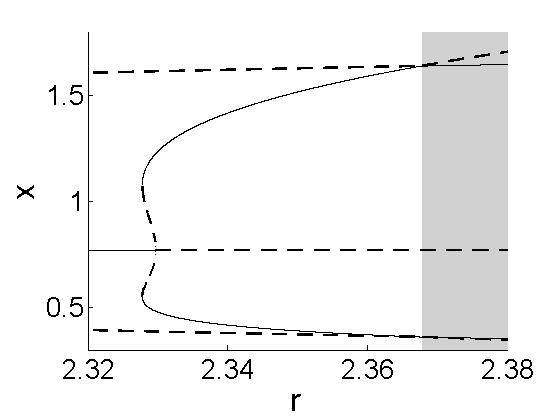}
     }
\hfill
\subfloat[$b=1.3$]{%
	\includegraphics[width =  0.48\linewidth]{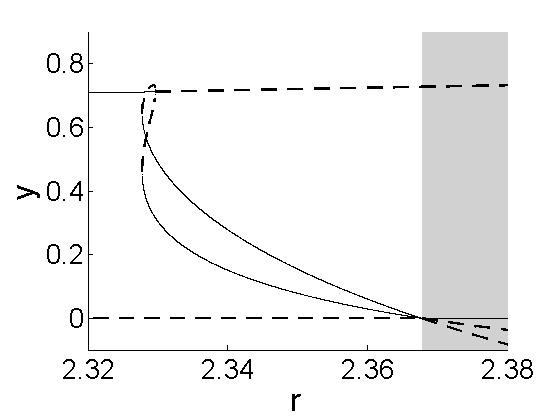}
	}
	\\
\subfloat[$b=1.5$]{%
       \includegraphics[width =  0.48\linewidth]{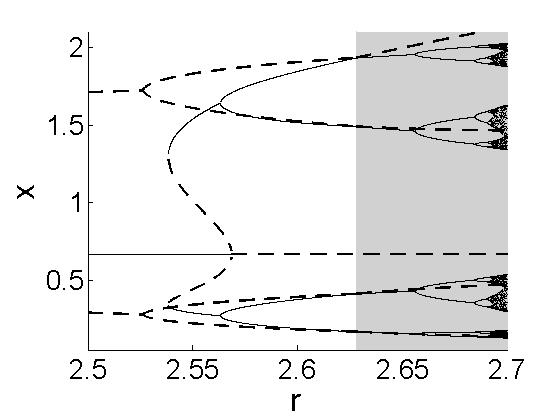}
     }
\hfill
\subfloat[$b=1.5$]{%
	\includegraphics[width =  0.48\linewidth]{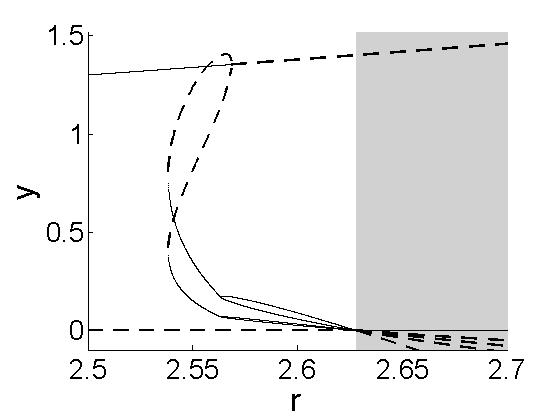}
	}
 \caption{Bifurcation diagrams for Model 2 for fixed $b$ as $r$ increases. Figures on the left show $x$ coordinates of stable (solid) and unstable (dashed) fixed points and cycles as $r$ increases. Figures on the right show $y$ coordinates of stable (solid) and unstable (dashed) fixed points and cycles as $r$ increases. Detailed descriptions of the dynamics and bifurcations are given in the text.}\label{Fig:Model2Bifurcations}
\end{figure}

For a fixed value of $b$ as $r$ increases, another attractor emerges. For certain values of $b$, there is a range of $r$ values for which bistability is observed. We describe the behavior for various fixed $b$ as $r$ is increased in the specified range noting that $R_0=\exp(r)$. This is illustrated in Figure \ref{Fig:Model2Bifurcations}.

\begin{itemize}
\item $b=1.004, 1.1, 1.2$ \\
For $r$ sufficiently high, the system has a stable interior equilibrium, an unstable equilibrium at $(1,0)$, and an unstable 2-cycle on the $x$-axis. As $r$ increases further, the interior equilibrium undergoes a supercritical flip bifurcation giving rise to a stable 2-cycle. As $r$ continues to increase, the 2-cycle moves towards the $x$-axis before colliding with the unstable 2-cycle and exchanging stability as it passes into the fourth quadrant.
\end{itemize}

%%%%%      FIGURE 4 HERE         %%%%%%%%%%%%%%%%%%%%%%%%%%%%%%%%%%%%%%%
%%%%%%%%%%%%%%%%%%%%%%%%%%%%%%%%%%%%%%%%%%%%%%%%%%
%%          FIGURE 4          %%%%%%%%%%%%%%%%%%%%%%%%%%%%%%%%%%%%%%%%%%%%%%%%%%%%%%%%%%%%%%%%%%%%%%%%%%%%%%%%%%%%%%

%%%%%%%%%%%%%%%%%%%%%%%%%%%%%%%%%%%%%%%%
%
%
%
\begin{itemize}
\item $b=1.3, 1.4$\\
For $r$ sufficiently high, the system has a stable interior equilibrium, an unstable equilibrium at $(1,0)$, and an unstable 2-cycle on the $x$-axis. As $r$ increases further, a stable 2-cycle emerges with an accompanying unstable 2-cycle in a saddle-node bifurcation of the second iterate of the mapping. Shortly thereafter, the interior equilibrium undergoes a subcritical flip bifurcation when the unstable 2-cycle in the first quadrant collides with it, and the equilibrium loses stability. For further discussion of subcritical flip bifurcations, see \cite{neubert1992subcritical} and \cite{whitley1983discrete}. As $r$ continues to increase, the stable 2-cycle moves towards the $x$-axis before colliding with the unstable 2-cycle on the $x$-axis and exchanging stability as it passes into the fourth quadrant.
\end{itemize}
%
%
%
%%%%%%%%%%%%%%%%%%%%%%%%%%%%%%%%%%%%%%%%%
%
\begin{itemize}
\item $b=1.5$\\
For $r$ sufficiently high, the system has a stable interior equilibrium, an unstable equilibrium at $(1,0)$, and an unstable 2-cycle on the $x$-axis. As $r$ increases further, we first observe that the unstable two-cycle on the $x$-axis period doubles into a four-cycle. Then, a stable 2-cycle emerges in the interior of the first quadrant with an accompanying unstable 2-cycle in a saddle-node bifurcation of the second iterate of the mapping. Then, the stable 2-cycle undergoes a period doubling bifurcation such that a stable 4-cycle emerges. Shortly thereafter, the interior equilibrium undergoes a subcritical flip bifurcation when the unstable 2-cycle in the interior of the quadrant collides with it. After this bifurcation, the coexistence equilibrium is unstable. As $r$ continues to increase, the stable 4-cycle moves towards the $x$-axis before colliding with the unstable 4-cycle and exchanging stability as it passes into the fourth quadrant.
\end{itemize}

\begin{figure}
\centering
\subfloat{%
       \includegraphics[width = 0.48\linewidth]{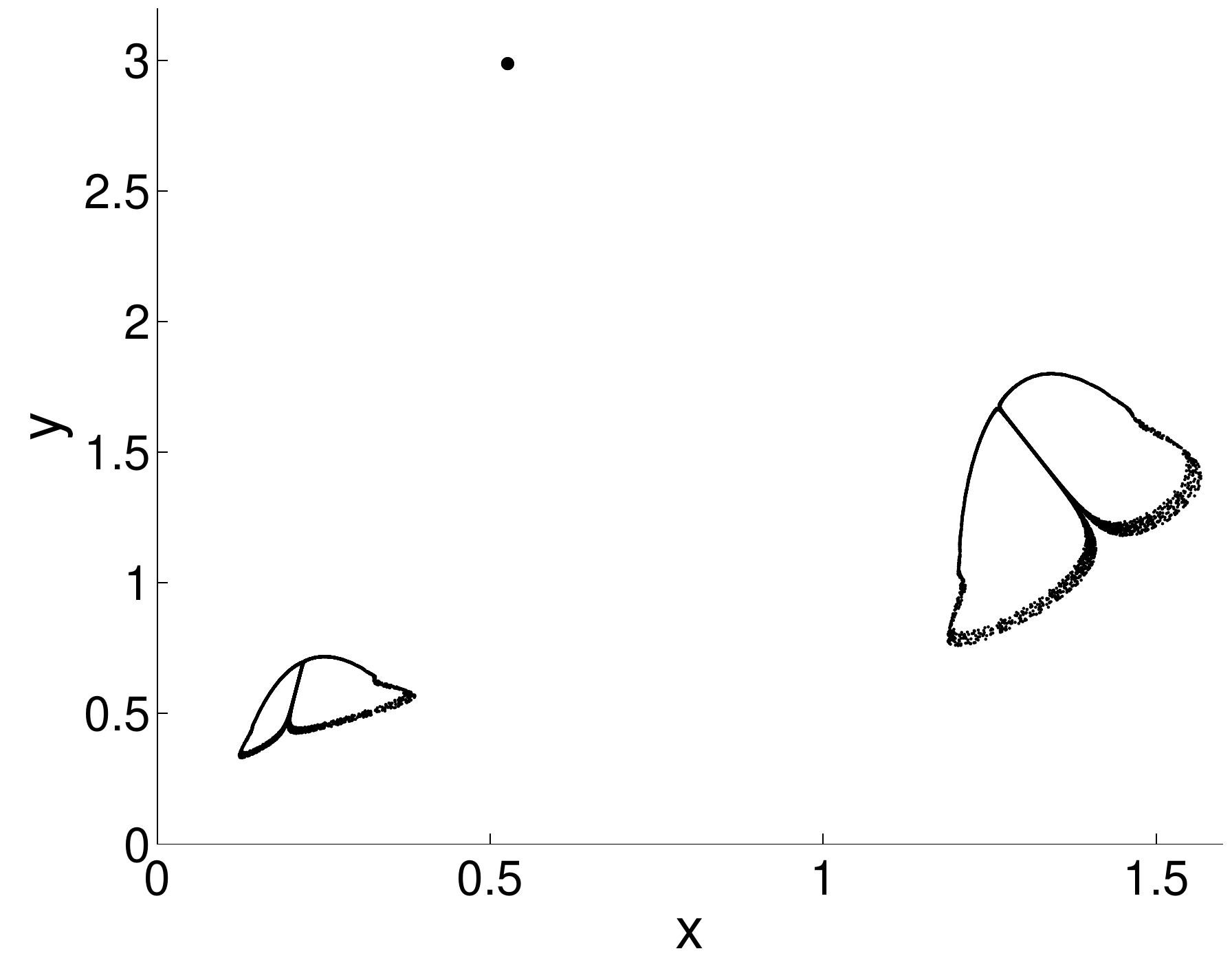}
     }
     \hfill
\subfloat{%
       \includegraphics[width = 0.48\linewidth]{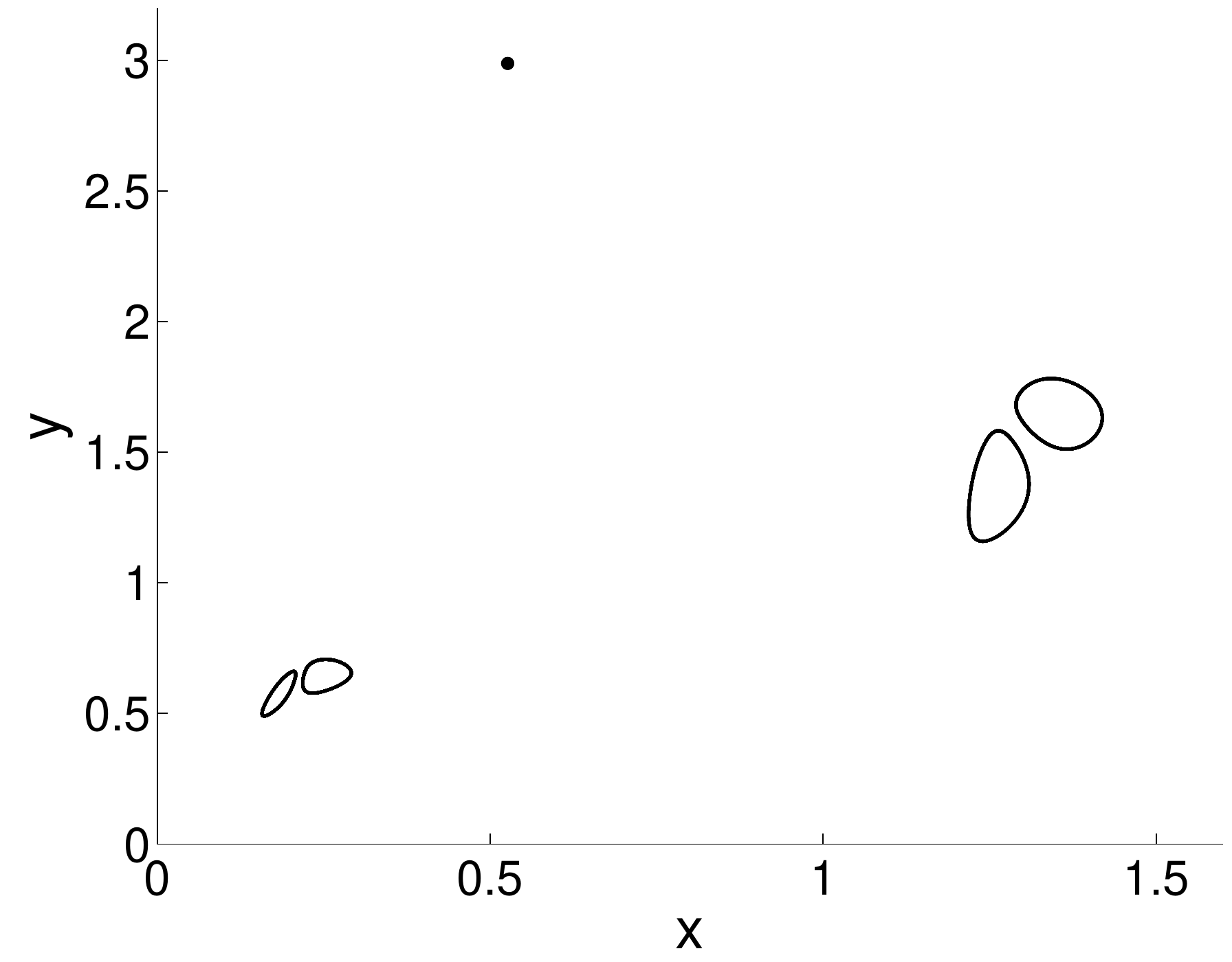}
     }
 \caption{Illustration of bistability between the equilibrium and another attractor for Model 2. The parameters for the left figure are $b=1.9,\ r=2.92$. For the right figure, $b=1.9,\ r= 2.9205$. For the figure on the left, the attractor is a region with fractal dimension. The attractor on the right consists of four circles such that the union of the circles is an invariant attracting set. If we continue to increase $r$ past $r=2.9205$, we see a 4-cycle that then period doubles to an 8-cycle. Initial conditions for both figures were $(0.8,0.7)$ for the equilibrium and $(0.3,0.4)$ for the other attractor. For clarity of the attractors, we ran $100,000$ iterations and plotted $30,000$ points for the attractors.}\label{fig:WeirdAttractor}
\end{figure}

It is evident that for higher values of $b$, the bifurcations associated with increasing $r$ are more complicated. Indeed, for $b= 1.9$,  $r=2.92$, the system has an attractor with fractal dimension. A small increase in $r$ to $r=2.9205$ results in an attractor made up of four circles such that the union of the four circles is an invariant attractor. In both cases, the equilibrium point is locally stable with its own basin of attraction. These two cases are shown in Figure \ref{fig:WeirdAttractor}. Further increases in $r$ result in a 4-cycle, which then period doubles into an 8-cycle.

%%%%%%           FIGURE 5 HERE        %%%%%%%%%%%%%%%%%%%%%%%%%
%%%%%%%%%%%%%%%%%%%%%%%%%%%%%%%%%%%%%%%%%%
%%%%%%%%%%        FIGURE 5        %%%%%%%%%%%%%%%%%%%%%%%%
%%%%%%%%%%%%%%%%%%%%%%%%%%%%%%%%%%%%

%%%%%%%%%%%%%%%%%%%%%%%%%%%%%%%%%%%%%%%%%%

%%%%%%%%%%%%%%%%%%%%%%%%%%%%%%%%%%%%%%%%%%%%%
\section{Model 3: Compensatory host density-dependence and exponential parasitism}\label{Section:M3}
For the third model under consideration, we return to fractional recruitment, equation (\ref{Eq:BevHoltxy}), and now incorporate a stronger parasitism term. That is, we now take the limit as $\kappa\to\infty$ in equation (\ref{Eq:NegBinom}), which results in exponential parasitism seen in equation (\ref{Eq:ExpPar}). Biologically, higher $\kappa$ corresponds to higher parasitoid aggregation, detailed in \cite{may1978host}.

The third model is
\begin{subequations}\label{Eq:Model3}
\begin{align}
x_{t+1}&=\left[\frac{R_0x_t}{1+(R_0-1)x_t}\right]e^{-y_t},\\
y_{t+1}&=b\left[\frac{R_0x_t}{1+(R_0-1)x_t}\right]\left(1-e^{-y_t}\right).
\end{align}
\end{subequations}
As with the other models, the coexistence equilibrium is in the interior of the first quadrant for $R_0>1$, $b>1$. This is shown in Section \ref{app:Model3FirstQuadrant}. For $b=1$, the coexistence equilibrium has collided with the exclusion equilibrium at $(1,0)$.  For $R_0=1$, system (\ref{Eq:Model3}) has a line of equilibria on the $x$-axis. As mentioned in Section \ref{Sec:CoexEq}, we cannot derive an explicit expression for the coexistence equilibrium for models with exponential parasitism.

\subsection{Stability region}
Even without an explicit expression for the coexistence equilibrium, we can determine the stability criteria. The first Jury condition is satisfied for $R_0>1$, $b>1$. Satisfying the first Jury condition is a sufficient condition for satisfying the second Jury condition. The third Jury condition, in turn, is satisfied in the $R_0$--$b$ plane below the curve
\begin{equation}
R_0=\frac{ye^{2y}}{e^y-1},\quad b=\frac{y^2e^{2y}-ye^y+y}{(e^y-1)(ye^y-e^y+1)},
\end{equation}
for positive $y$. We determined this parametric curve for the third Jury condition, inequality (\ref{Eq:Jury3}), by solving the three equations
\begin{align}
\Delta &= 1,\\
u(x,y) = g(x)[1-yh(y)]&=1,\label{u1}\\
v(x,y) = bxg(x)h(y)&=1,\label{v1}
\end{align}
to eliminate $x$ and write $b$ and $R_0$ as functions of $y$. Equations (\ref{u1}) and (\ref{v1}) are the equations for the host and parasitoid nullclines given in system (\ref{Sys:uv1}). Details for all three Jury conditions are given in Sections \ref{app:Model3Jury1}--\ref{app:Model3Jury3}.

As shown in the stability region in Figure \ref{fig:Model3StabilityRegion}, the Jury 3 condition is the interesting feature of the stability region for Model 3. Crossing the dashed curve in parameter space from below corresponds to violating the third Jury condition such that both eigenvalues leave the unit disc in the complex plane. For our model, this yields a supercritical Neimark--Sacker bifurcation \cite{wiggins2010introduction}, where the equilibrium loses stability and is replaced by a stable, quasiperiodic attractor that is topologically similar to a circle. These attractors are commonly referred to as invariant circles. The bifurcation itself is sometimes referred to as a discrete Hopf bifurcation. Crossing the line $b=1$ violates the first Jury condition and results in a transcritical bifurcation. Crossing the line $R_0=1$ violates both the first and third Jury conditions.
%

%%%%%%%%   FIGURE 6 HERE  %%%%%%%%%%%%%%%%%%%%%%
%%%%%%%%%%%%%%%%%%%%%%%%%%%%
%%           FIGURE 6       %%%%%%%%
%%%%%%%%%%%%%%%%%%%%%%%%%%%%%%%%%%
\begin{figure}
\centering
\subfloat{%
       \includegraphics[width = 0.48\linewidth]{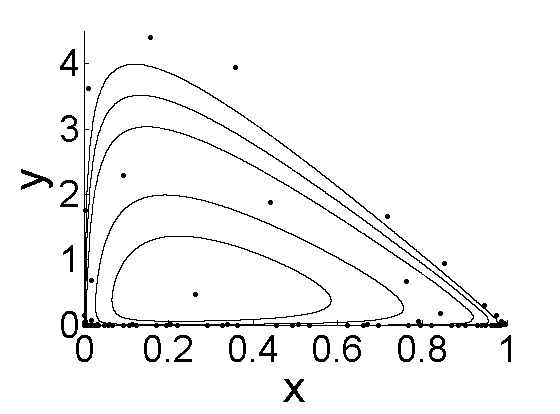}
     }
\hfill
\subfloat{%
       \includegraphics[width = 0.48\linewidth]{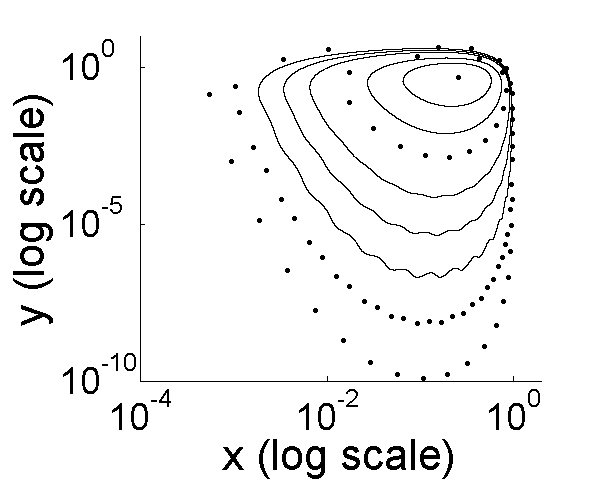}
     }
 \caption{To investigate the interior attractor for Model 3, we fix $R_0=2$ and increase $b$. The value of $b$ that produces the stable equilibrium in these figures is $b=3$. Moving outward from this equilibrium point, the invariant circles and phase-locked cycles correspond to $b=3.5, 4, 4.41, 5, 5.5, 6, 6.5, 7$. Crossing the boundary of the stability region results in a Neimark--Sacker bifurcation such that an invariant circle becomes the stable attractor. The invariant circle grows and undergoes phase-locking alternating with invariant circles as $b$ continues to increase. As $b$ increases, the lower portion of the attractor approaches the $x$-axis.}\label{fig:Model3Bif}
\end{figure}

\subsection{Bifurcations and attractors}\label{sec:M3bifn}
In order to illustrate the bifurcations and types of attractors for different parameter choices, we fix $R_0=2$ and increase $b$. Figure \ref{fig:Model3Bif} shows the attractors for increasing $b$ values. For $b$ values below the Jury 3 curve in Figure \ref{fig:Model3StabilityRegion}, the coexistence  point is stable. After the Neimark--Sacker bifurcation, the complex eigenvalues of the fixed point are larger than one in magnitude, and the attractor is either a quasiperiodic invariant circle or a periodic $n$-cycle, increasing in size as $b$ increases. 

We now qualitatively describe the behavior of the system for parameters outside the stability region. For some values of $b$ just above the Jury 3 curve, the system has a stable invariant circle. For other values of $b$ also just above the Jury 3 curve, there is a pair of periodic orbits on the invariant circle, one stable and one unstable. When the rotation number of the periodic orbits is $p/q$, the system has a $p/q $ resonance \cite{aronson1982bifurcations}. Specifically, we consider the case of weak resonance such that $q\neq 1,2,3,4$ as the eigenvalues pass through the unit circle \cite{whitley1983discrete}. The set of parameter values for which the system has a periodic orbit with rational rotation number $p/q$ is known as an Arnold tongue \cite{wiggins2010introduction}. 

%%%%%%%%%%          FIGURES 7 & 8 HERE     %%%%%%%%%%%%%%%%

%%%%%%%%%%%%%%%%%%%%%%%%%%%%%%%
%%%%%          FIGURE 7        %%%%%%%%%%%%%%
%%%%%%%%%%%%%%%%%%%%%%%%%%%%%%%%%%%%%%%%%%%%%%%%%

\begin{figure}
\centering
\includegraphics[width=0.98\linewidth]{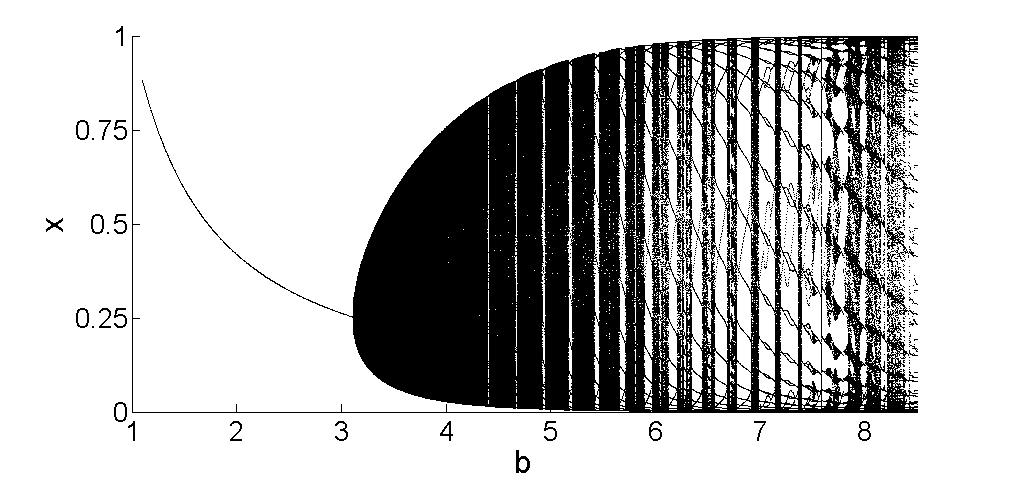}\\
\caption{Model 3 bifurcation diagram illustrating the $x$ coordinates of the stable attractor for $R_0=2$ and varying $b$. Note the Neimark--Sacker bifurcation that occurs when the equilibrium loses stability and an invariant circle becomes the attractor, corresponding to multiple $x$ values for a single value of $b$. (Bifurcation diagram for $y$ not shown here.)}\label{fig:3xBif} 
\end{figure}
%%%%%%%%%%%%%%%%%%%%%%%%%%%
%%%%%%%%%%%%%%%%%%%%%%%%%%%%%%%%%%%%%
%%%%%%          FIGURE 8            %%%%%%%%%%%%%%%%%%%%%%
%%%%%%%%%%%%%%%%%%%%%%%%%%%%%%%%%%%%%%%
\begin{figure}
\centering
\includegraphics[width=0.98\linewidth]{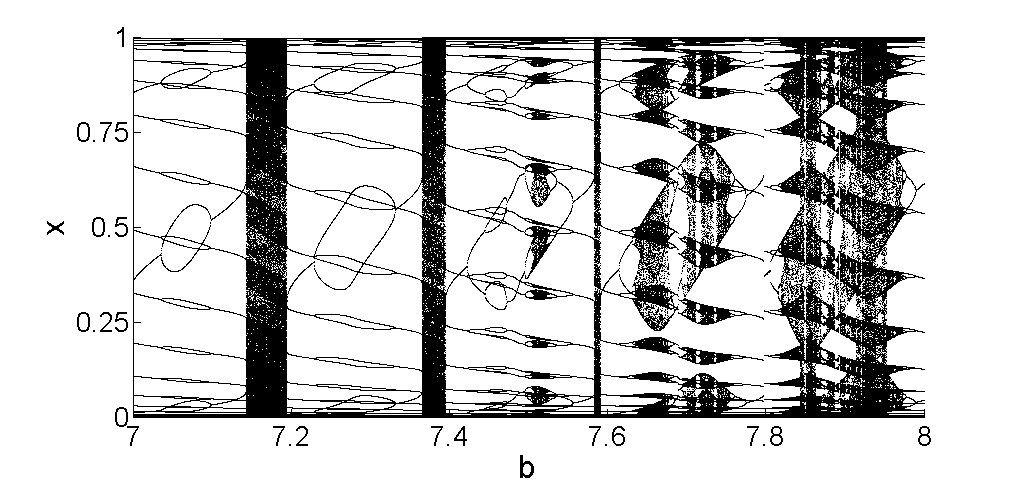}\\
\caption{Model 3 bifurcation diagram with a narrower range of $b$ values. As $b$ increases, phase locking occurs along with period doubling and halving, interspersed with regions of invariant circles corresponding to a dense set of $x$ coordinates of the attractor. These phenomena occur as the eigenvalues pass through the resonance horns corresponding to phase-locked $n$-cycles.}\label{fig:3xBifDetail}
\end{figure}
%%%%%%%%%%%%%%%%%%%%%%%%%%%%%%

Alternately, Arnold tongues or resonance horns may describe a cusped region in the complex plane where eigenvalues within the horn correspond to the existence of a stable periodic orbit with rational rotation number \cite{aronson1982bifurcations, lauwerier1986two, kuznetsov2004elements}. The eigenvalue will typically intersect an infinite number of these resonance horns near the unit circle \cite{whitley1983discrete}. In our case, as $b$ continues to increase, the eigenvalues of the coexistence point pass in and out of resonance horns or Arnold tongues. Whenever the eigenvalues are within a resonance horn, the system is phase-locked, and we observe a stable $n$-cycle in the $x$-$y$ plane. As the eigenvalues continue to grow in magnitude, the Arnold tongues are wider, and there are broader windows of phase-locking in the bifurcation diagram as the parameter $b$ increases. Within these windows, the system may undergo changes to the period of the $n$-cycle as the eigenvalues enter and exit overlapping resonance horns with differing rational rotation numbers. This behavior is visible in  Figures \ref{fig:3xBif} and \ref{fig:3xBifDetail}, which show bifurcation diagrams of the $x$ coordinates of the attractors to illustrate the changes in the system as $b$ increases for fixed $R_0=2$.  

Returning to Figure \ref{fig:Model3Bif}, we note that as $b$ increases, the stable attractor in the system grows, and the lower portion approaches the $x$-axis. While the interior of the first quadrant is invariant for system (\ref{Eq:Model3}), numerical simulations of the system for values of $b$ much past $b=8.5$ result in rounding small positive values of $y$ down to identically 0. Thus, numerical simulations are limited in their ability to demonstrate the behavior of the system for even finite parameter values. It is ecologically likely, however, that for sufficiently small values of $y$, stochastic events would wipe out the parasitoid population, after which, the dynamics of the system would reflect the dynamics observed on the $x$-axis.
%%%%%%%%%%%%%%%%%%%%%%%%%%%%%%%%%%%%%%%
\section{Model 4: Overcompensatory host density-dependence and exponential parasitism}\label{Section:M4}
The fourth model uses an exponential form for both density dependence and parasitism. This corresponds to equations (\ref{Eq:Ricker}) and (\ref{Eq:ExpPar}). The model is thus
\begin{subequations}\label{Eq:Model4}
\begin{align}
x_{t+1}&=x_t e^{r(1-x_t)}e^{-y_t},\\
y_{t+1} &= b x_t e^{r(1-x_t)}\left(1-e^{-y_t}\right) .
\end{align}
\end{subequations}
Because of the exponential parasitism term, there is not an explicit expression for the coexistence equilibria solutions to system (\ref{Eq:Model4}).

Unlike the previous models, $b>1$ is not necessary for the occurrence of a coexistence equilibrium in the interior of the first quadrant. As shown in Figure \ref{fig:Model4StabilityRegion}, there is a region above $r=2$ and below $b=1$ for which there are two coexistence equilibria, only one of which may be stable.  For $0<r<2$, when $b=1$, the single coexistence equilibrium has collided with the exclusion equilibrium at $(1,0)$. For $r=0$ ($R_0=1$), system (\ref{Eq:Model4}) has a line of equilibria on the $x$-axis. 

This model was previously studied by Kang et al. \cite{kang2008dynamics} in the context of a plant--herbivore system. However, our nondimensionalization and methods of analysis differ from theirs. In particular, Kang et al. \cite{kang2008dynamics} studied stability of the equilibria numerically, while we use analytic methods. This allows us to find a bifurcation that is missing from their analysis, discussed below.

%%%%%%%%%%%%%%%%%%%%%%%%%%%%%%%%%%%%%%%%%%%%%%%
\subsection{Stability Region}\label{Section:M41}
This model uses exponential forms for both density-dependent recruitment and parasitism. The stronger nonlinearity in density dependence and the stronger form of parasitism result in both the second and third Jury conditions functioning as interesting boundaries of the stability region, seen in Figure \ref{fig:Model4StabilityRegion}. The stability conditions are:
\begin{enumerate}
\item For $0<r<2$, Jury condition 1 is satisfied above $b=1$; for $r>2$, Jury condition 1 is satisfied above the curve
\begin{equation}\label{Eq:M4Jury1}
r=\frac{y^2e^y}{1+ye^y-e^y},\quad b= \frac{y^2e^y}{\left(e^y-1\right)^2},
\end{equation}
\item Jury condition 2 is satisfied above the curve
\begin{equation}\label{Eq:M4Jury2}
r=\frac{e^y(y^2+2y+2)-2}{e^y(y+1)-1},\quad b=\frac{2y(e^y-1)+y^2e^y(2+y)}{(2+y)e^{2y}-4e^y-y+2},
\end{equation}
\item Jury condition 3 is satisfied for $r>0$ ($R_0>1$) and below the curve
\begin{equation}\label{Eq:M4Jury3}
r=\frac{e^y(y^2+y-1)+1}{ye^y}, \quad b=\frac{e^y(y^3+y^2-y)+y}{\left(e^y-1\right)\left(ye^y-e^y+1\right)}.
\end{equation}
\end{enumerate}
The parametric curves (\ref{Eq:M4Jury1}), (\ref{Eq:M4Jury2}) and (\ref{Eq:M4Jury3}) are all defined for positive $y$. When $b=1$ and $0<r<2$, the first Jury condition is violated. When $r=0$ ($R_0=1$), the first and third Jury conditions are violated. 

Details for determining all three conditions are given in Sections \ref{app:Model4Jury1}--\ref{app:Model4Jury3}. For each parametrically-defined curve, we used the equations for the host and parasitoid nullclines, equations (\ref{Eq:uis1}) and (\ref{Eq:vis1}), with either $1-\tau+\Delta = 0$, $1+\tau+\Delta =0$, or  $\Delta =1$ to eliminate $x$ and write $b$ and $r$ as functions of $y$.

Note that curve (\ref{Eq:M4Jury1}) is entirely below curve (\ref{Eq:M4Jury2}) (not shown). Curves (\ref{Eq:M4Jury1}) and (\ref{Eq:M4Jury2})  are visibly indistinguishable at the scale used in Figure \ref{fig:Model4StabilityRegion}. For a given $r>2$, the value of $b$ must be above curve (\ref{Eq:M4Jury2}) for stability to be guaranteed. The dynamics of the system for parameters between curves (\ref{Eq:M4Jury1}) and (\ref{Eq:M4Jury2}) are discussed in Section \ref{Section:M3Bifn_etc} below.

Returning to Figure \ref{fig:Model4StabilityRegion}, we consider the bifurcations that occur when the Jury conditions are violated. Crossing the dotted curve from above violates the second Jury condition, and the system undergoes a subcritical period-doubling bifurcation. Crossing the dashed curve from below corresponds to a supercritical Neimark--Sacker bifurcation \cite{wiggins2010introduction}, where the equilibrium loses stability and is replaced by a stable, invariant circle. Crossing the solid horizontal line, $b=1$, $0<r<2$ from above corresponds to a transcritical bifurcation where the unique coexistence equilibrium collides with the exclusion equilibrium on the $x$-axis.

While Kang et al. \cite{kang2008dynamics} use a different nondimensionalization, their parameter $a$ is the same as our parameter $b$, the product of searching efficiency, parasitoid clutch size, and the host carrying-capacity. Thus, it holds from their work that for $b>1$, system (\ref{Eq:Model4}) has a unique coexistence equilibrium. For $r>2$, our results indicate that above the first Jury condition curve, (\ref{Eq:M4Jury1}), and below $b=1$, there are two coexistence equilibria. For the region above the second Jury condition curve, (\ref{Eq:M4Jury2}), below both $b=1$ and below the third Jury condition curve, (\ref{Eq:M4Jury3}), the equilibrium point with the larger $y$ value is stable. This region extends infinitely in the $r$ direction since the third Jury condition curve, (\ref{Eq:M4Jury3}), remains above the second Jury condition curve,  (\ref{Eq:M4Jury2}), even after the curve (\ref{Eq:M4Jury3}) is below $b=1$. See Figure \ref{fig:Model4StabilityRegion} and details in Sections \ref{app:Model4Jury1}--\ref{app:Model4Jury3}. The second Jury condition curve,  (\ref{Eq:M4Jury2}), is the one that was missed by Kang et al. \cite{kang2008dynamics}. 

%%%%%%%%%%%%%%%%%%%%%%%%%%%%%%%%%%%%%%%%%%%%
\subsection{Bifurcations and Attractors}
\label{Section:M3Bifn_etc}

%%%%%%%%%%%%%%%%%%%%%%%%
%%%%%%%%        FIGURES  9 & 10 & 11 HERE        %%%%%%%%%%%%%%%%%%%%%%%%%%%%%%%%
%%%%%%%%%%%%%%%%%%%%%%%%%%%%%%%%%
%%%%%%              FIGURE 9                %%%%%%%%%%%%%
%%%%%%%%%%%%%%%%%%%%%%%%%%%%%
\begin{figure}
\centering
\subfloat[$b=0.959$]{%
       \includegraphics[width =  0.32\linewidth]{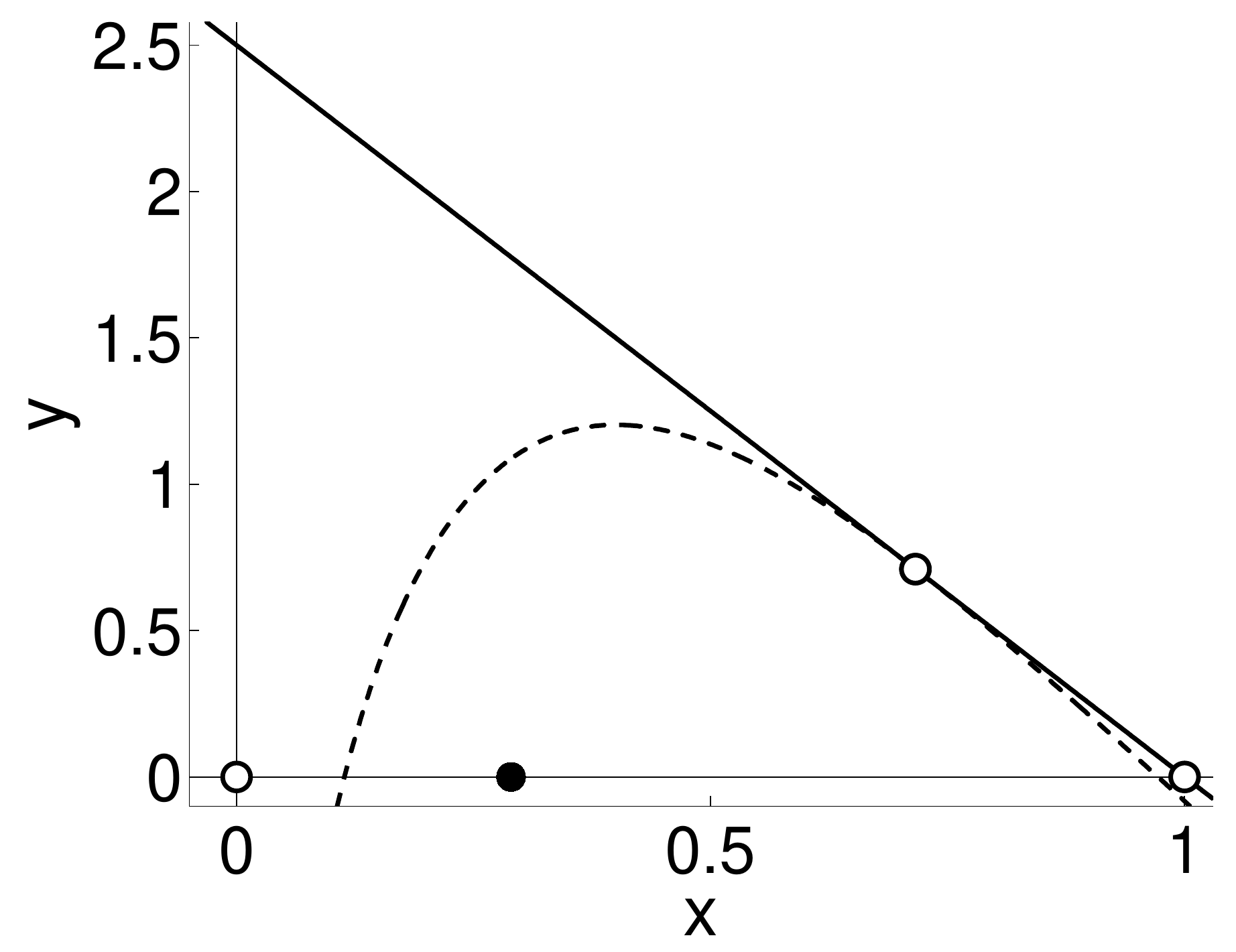}\label{fig:Tangency}
     }
\hfill
\subfloat[$b=0.96$]{%
       \includegraphics[width = 0.32\linewidth]{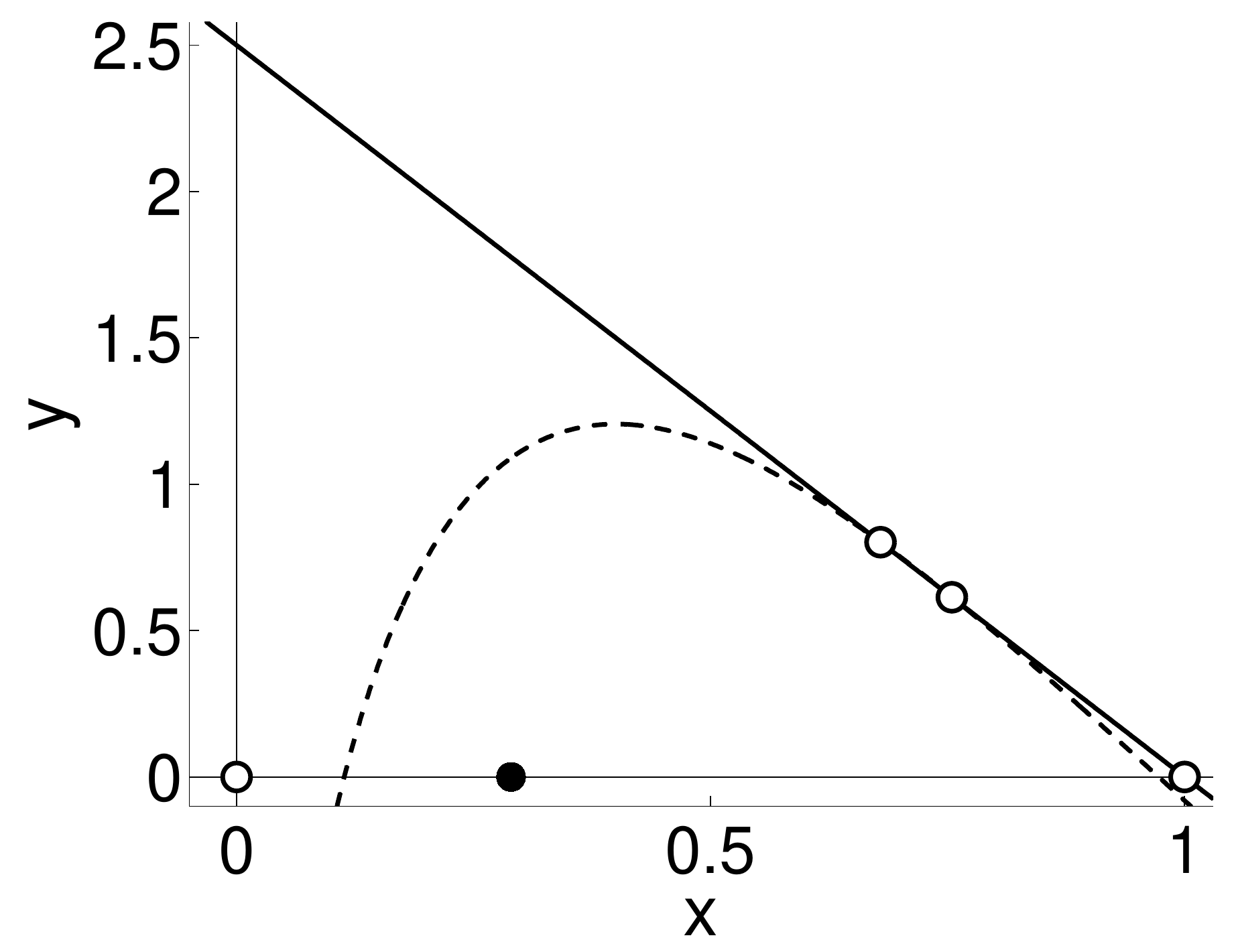}
     }
\hfill
\subfloat[$b=0.9615$]{%
       \includegraphics[width =  0.32\linewidth]{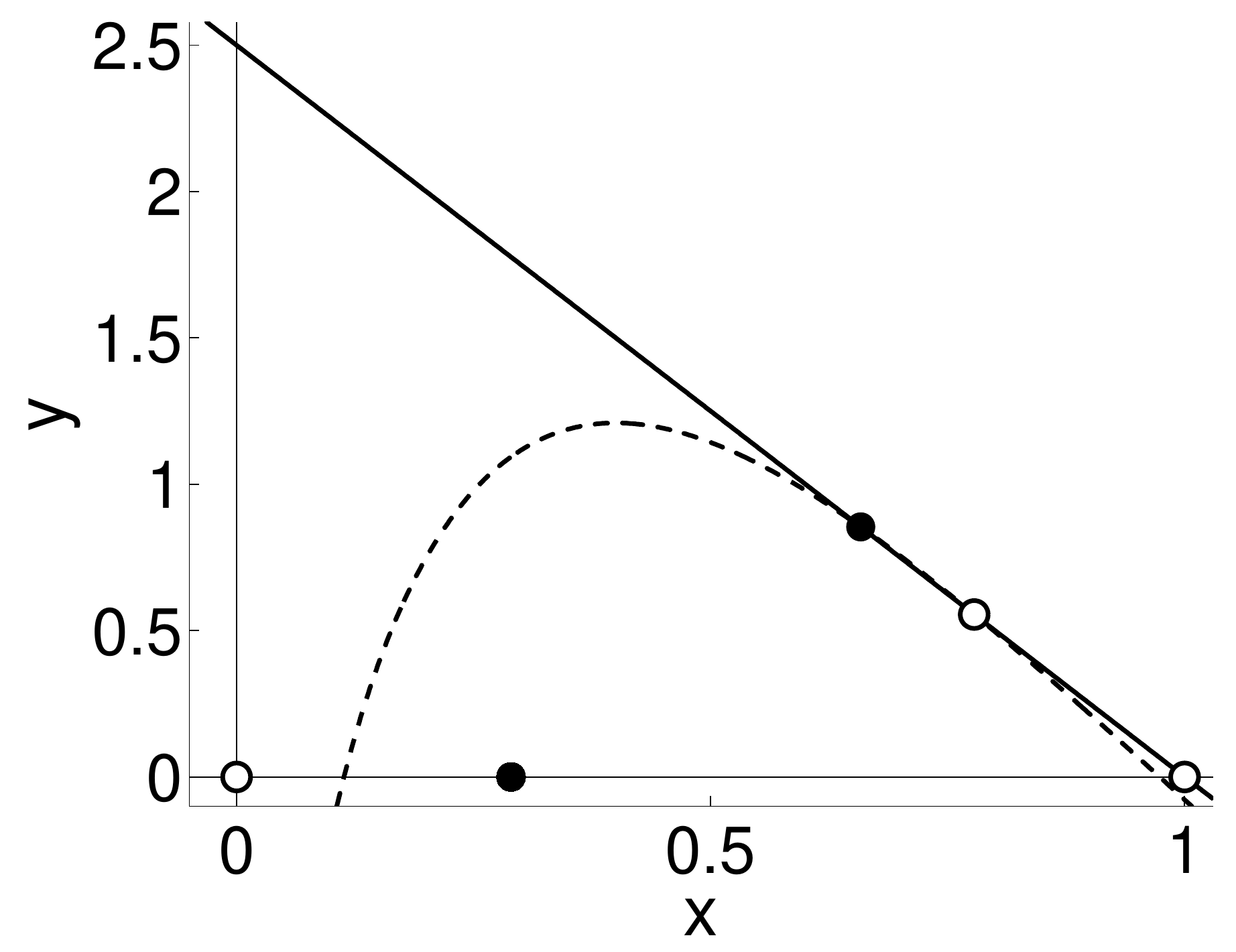}
     }
	\\
\subfloat[$b=0.98$]{%
	\includegraphics[width =  0.32\linewidth]{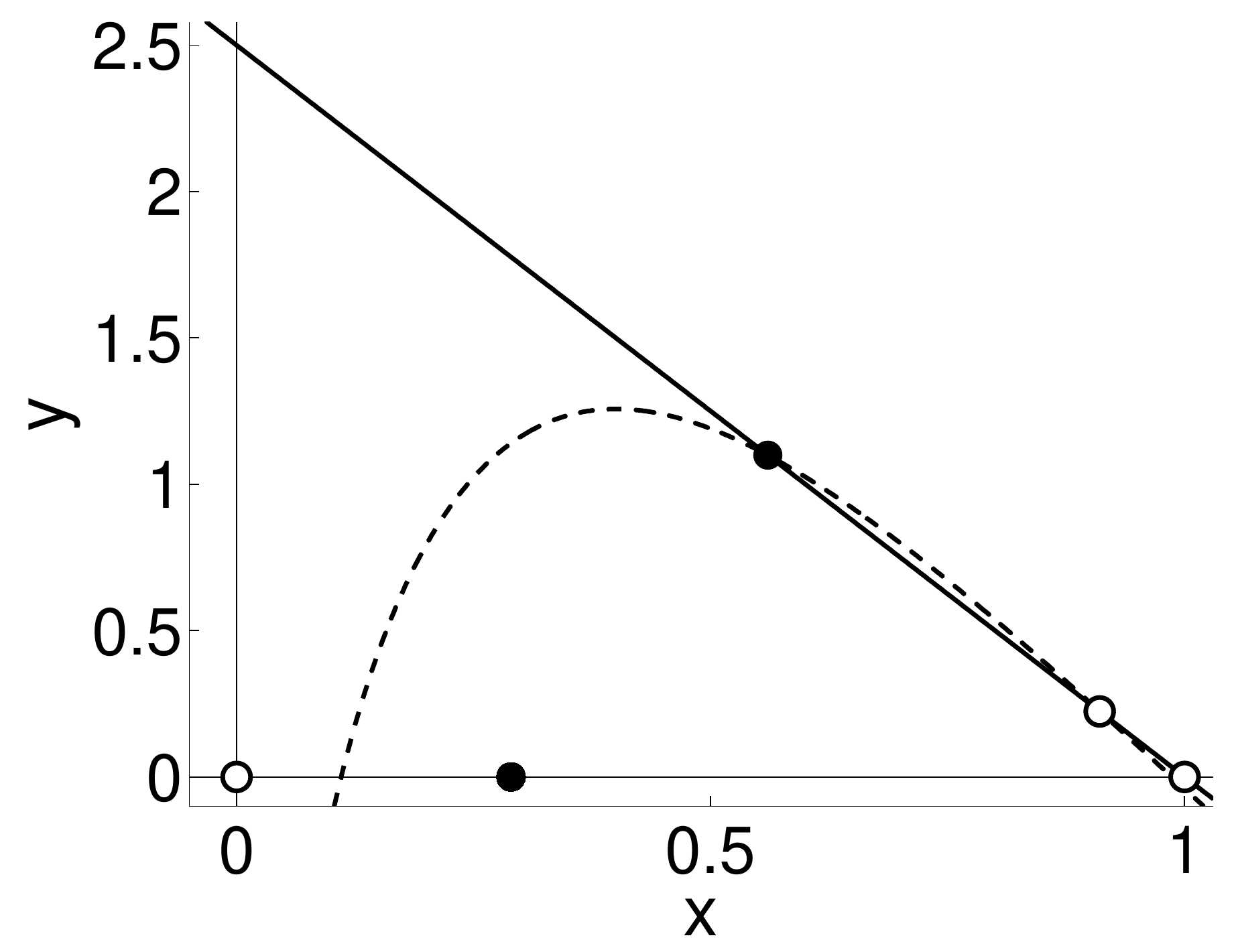}\label{fig:AlmostCollision}
	}
	\hfill
\subfloat[$b=1.05$]{%
       \includegraphics[width =  0.32\linewidth]{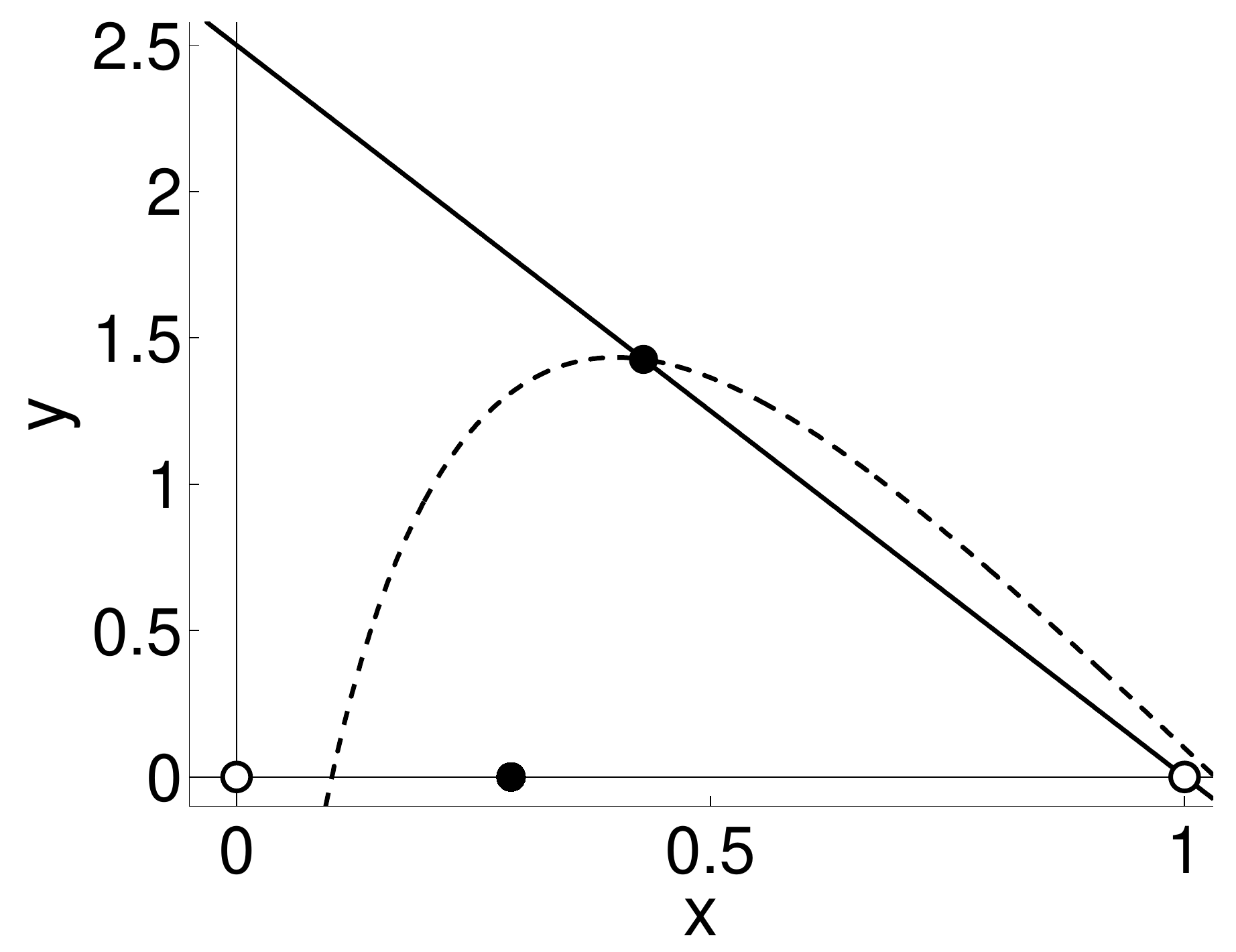}
     }
\hfill
\subfloat[$b=1.6$]{%
	\includegraphics[width =  0.32\linewidth]{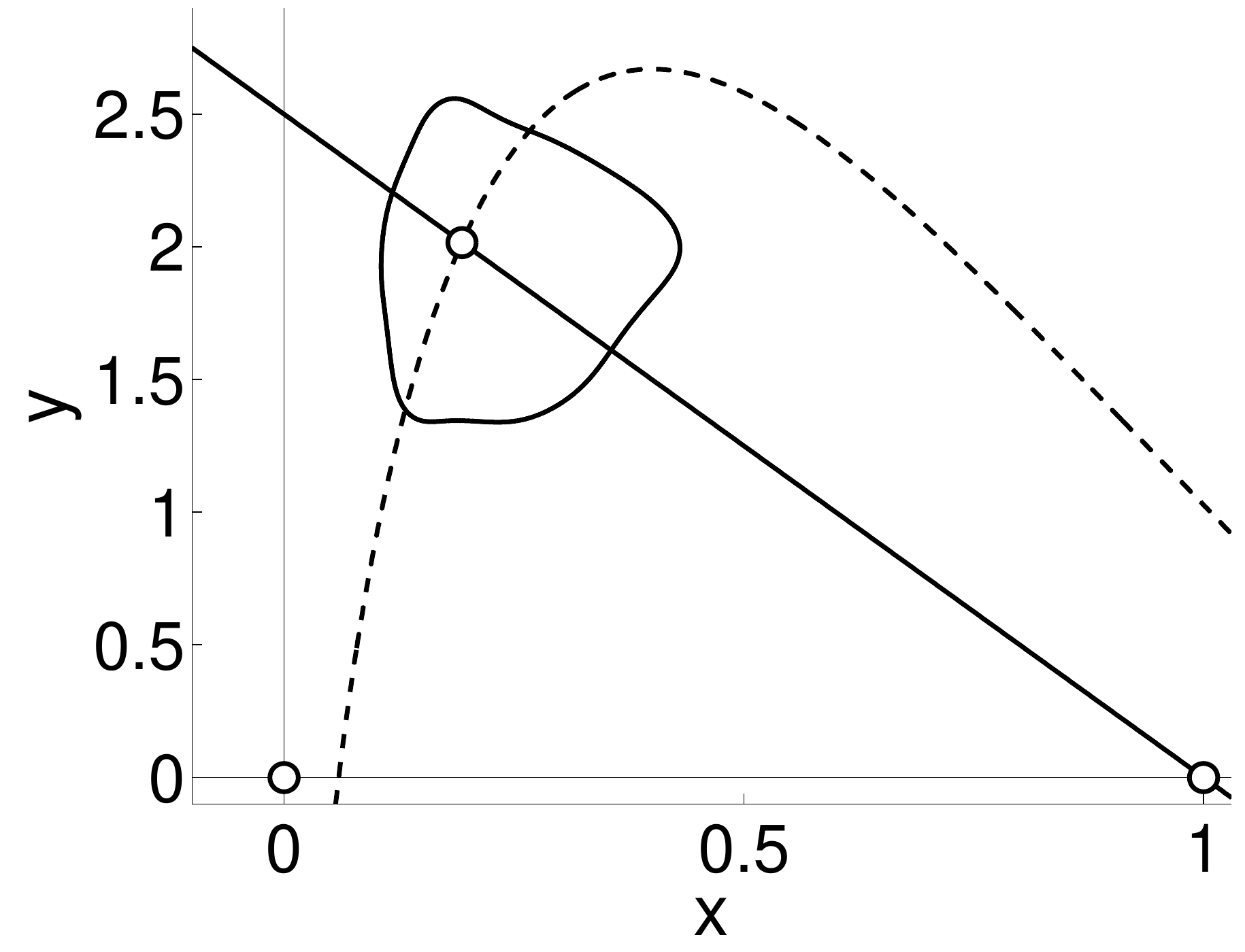}
	}
 \caption{Phase portraits for $r=2.5$ ($R_0\approx 12.18$) as $b$ increases. The host nullcline is shown with the solid line. As $b$ increases, the parasitoid nullcline, shown with a dotted line, changes shape. Attractors are shown with filled circles while unstable equilibria are shown with open circles. At $b\approx 0.959$, a saddle-node bifurcation occurs. Note the tangency between the nullclines. This condition is equivalent to the condition for the first Jury condition to be satisfied. Both coexistence equilibria are initially unstable, but after the subcritical period-doubling bifurcation seen in Figure \ref{fig:4yBifZoom}, the upper of the two equilibria is stable. The lower coexistence equilibrium moves towards the $x$-axis and collides with the exclusion equilibrium at $b=1$. For $b>1$, the coexistence equilibrium is unique. When the third Jury condition is violated, the coexistence equilibrium undergoes a supercritical Neimark--Sacker bifurcation. A stable quasiperiodic invariant circle becomes the attractor. Note that the solid circle on the $x$-axis is one of two points that form a stable two-cycle on the axis.}\label{Fig:Model4xyBifurcations}
\end{figure}
%%%%%%%%%%%%%%%%%%%%%%%%%%%%%%%%%
%%%%%%%%%%        FIGURE 10       %%%%%%%%%%%%%
%%%%%%%%%%%%%%%%%%%%%%%%%%%%%%
\begin{figure}
\centering
\includegraphics[width=0.98\linewidth]{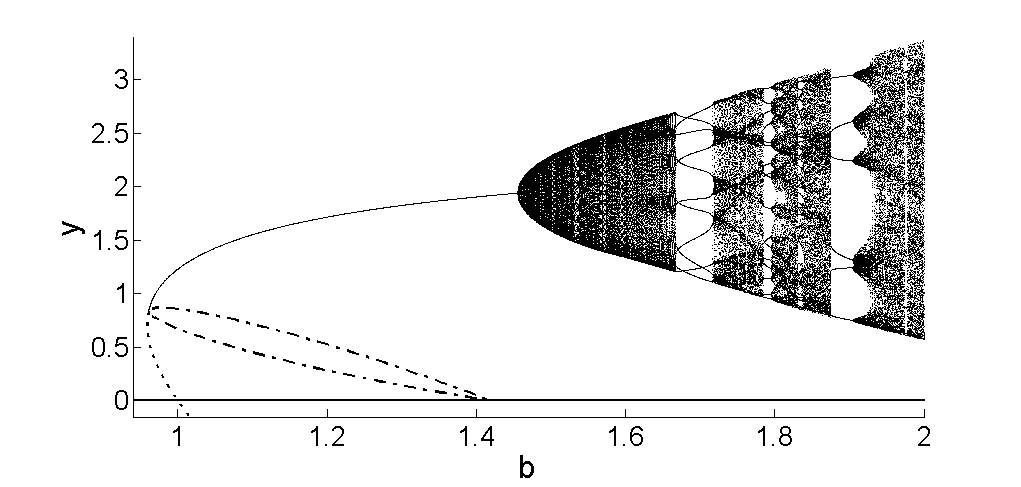}\\
\caption{Model 4 bifurcation diagram illustrating the $y$ coordinates of the stable attractor for $r=2.5$ ($R_0\approx 12.18$) and varying $b$. The fixed points emerge in a saddle-node bifurcation as $b$ crosses the first Jury condition curve (\ref{Eq:M4Jury1}). The upper of the two equilibria undergoes a subcritical period-doubling bifurcation in which it becomes stable. The resulting unstable two cycle is shown with the dash-dot line. The dotted line is the unstable equilibrium, which crashes through exclusion equilibrium on the $x$-axis at $b=1$. The unstable two-cycle also crashes through the $x$-axis. The stable equilibrium loses stability through a Neimark--Sacker bifurcation, and an invariant circle becomes the attractor, corresponding to multiple $y$ values for a single value of $b$. (Bifurcation diagram for $x$ not shown here.)}\label{fig:4yBif}
\end{figure}
%%%%%%%%%%%%%%%%%%%%%%%%%%%%%%%%%%%%%%%

%%%%%%%%%%%%%%%%%%%%%%%%%%%%%%%%%%%
%%%%%%%%                FIGURE 11          %%%%%%%%%%%%%%%%
%%%%%%%%%%%%%%%%%%%%%%%%%%%%%%%%
\begin{figure}
\centering
\subfloat{%
 \includegraphics[width = 0.48\linewidth]{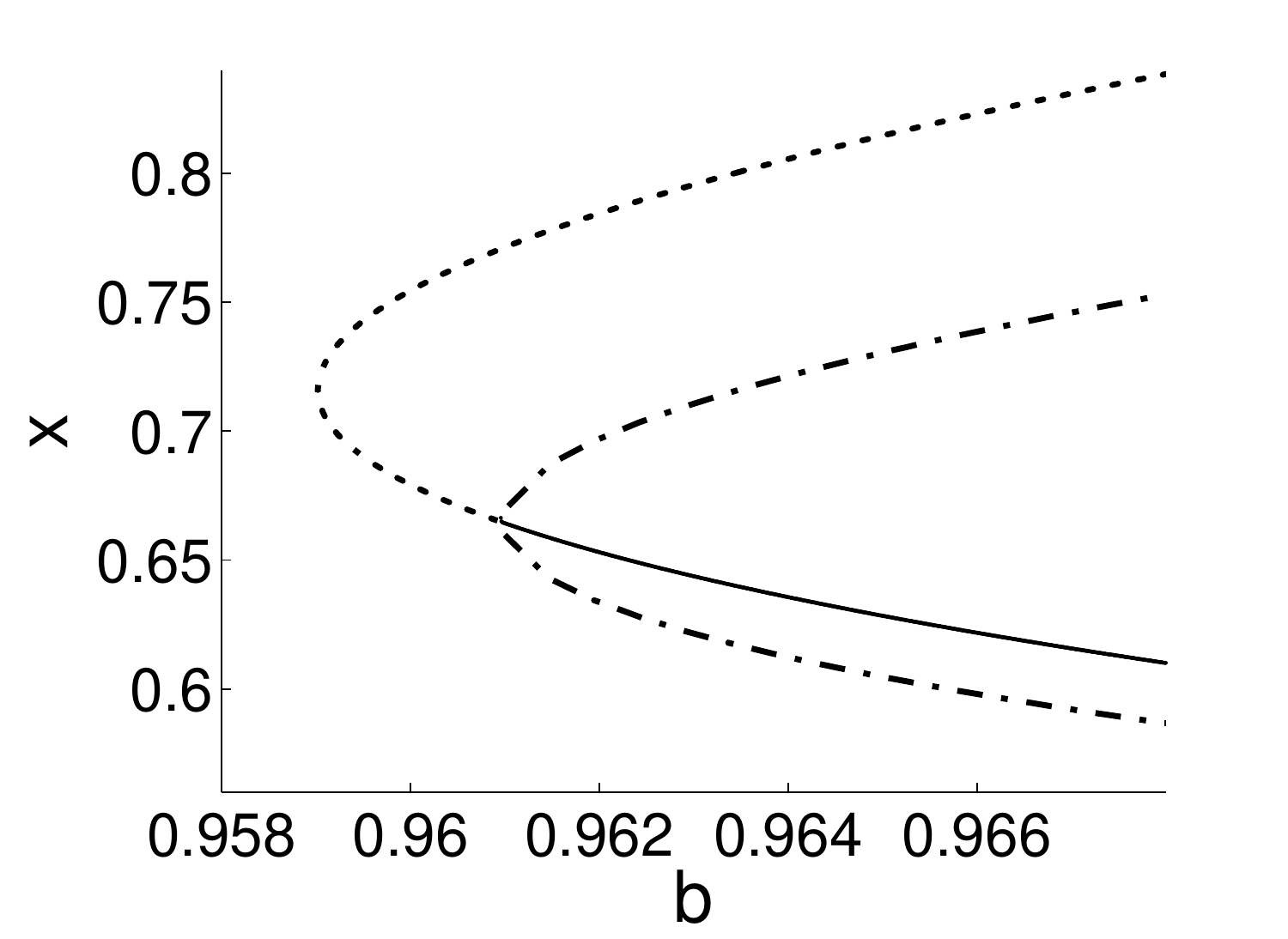}
     }
\hfill
\subfloat{%
\includegraphics[width=0.48\linewidth]{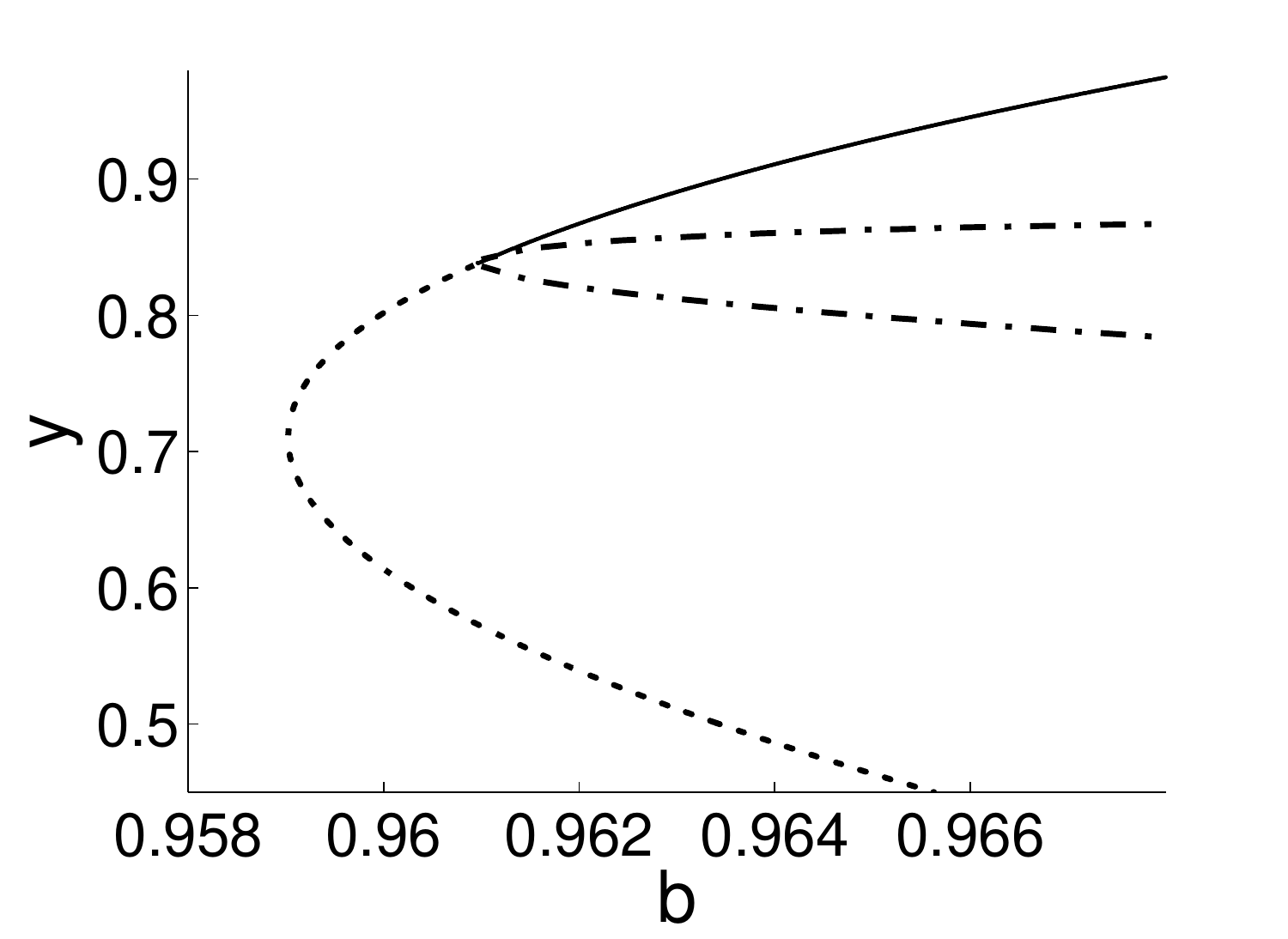}
     }
\caption{Model 4 bifurcation diagram with a much narrower range of $b$ values, again for $r=2.5$ ($R_0\approx 12.18$). Here, the saddle-node bifurcation is clearly visible such that both equilibria are initially unstable.  The left figure shows the $x$ coordinates and the right figure shows the $y$ coordinates for the same range of $b$ values. Dotted lines correspond to unstable equilibria.The equilibrium with the larger $y$ value undergoes a subcritical period-doubling bifurcation and gains stability as an unstable two cycle is born, shown with a dash-dot line. This behavior was missed in the numerical investigations by Kang et al. \cite{kang2008dynamics} but can be found analytically from equations (\ref{Eq:M4Jury1}) and (\ref{Eq:M4Jury2}). The solid line corresponds to where the equilibrium with the larger $y$ value is stable.}\label{fig:4yBifZoom}
\end{figure}
%%%%%%%%%%%%%%%%%%%%%%%%%%

In order to clearly illustrate the bifurcations and dynamics of the system, we fix $r=2.5$ ($R_0\approx 12.18$) and increase $b$. We have chosen a value of $r$ for which the second Jury condition curve is the lower boundary of the stability region, seen in Figure \ref{fig:Model4StabilityRegion}. The host and parasitoid nullclines are shown in Figure \ref{Fig:Model4xyBifurcations} with stable and unstable equilibria and other attractors for selected values of $b$. A bifurcation diagram for increasing $b$ values is shown in Figure \ref{fig:4yBif}. 

For $r=2.5$ and $b$ below the first Jury condition curve, equation (\ref{Eq:M4Jury1}), there are no coexistence equilibria. When we increase $b$ to the first Jury condition curve, the host and parasitoid nullclines are tangent, seen in Figure \ref{fig:Tangency}. For slightly higher values of $b$, both coexistence equilibria are unstable. This differs from the claim made in Kang et al. \cite{kang2008dynamics} that one equilibrium is stable after the saddle-node bifurcation. However, the instability of both coexistence equilibria occurs for a tiny range of $b$ values, from $0.959<b<0.961$. The upper of the two equilibria undergoes a subcritical period-doubling bifurcation and gains stability as $b$ crosses the dotted curve shown in Figure \ref{fig:Model4StabilityRegion}, the second Jury condition curve. The resulting unstable two-cycle was found numerically and is shown in Figures \ref{fig:4yBif} and \ref{fig:4yBifZoom}.

We continue with the bifurcations as $b$ increases past $b=1$. Returning to Figure \ref{Fig:Model4xyBifurcations}, we see that at $b=1$, the lower of the coexistence equilibria collides with the exclusion equilibrium as it passes into the fourth quadrant. As $b$ continues to increase, the unstable two-cycle in the interior of the first quadrant eventually crashes through the $x$-axis, passing into the fourth quadrant.  For sufficiently low values of $b$, the two-cycle on the axis is a competing stable attractor. When $b$ crosses the third Jury condition curve, a Neimark--Sacker bifurcation results and a quasiperiodic stable invariant circle is born. As seen in Figure \ref{fig:4yBif}, the complex eigenvalues of the coexistence equilibrium point again pass in and out of Arnold tongues, resulting in phase-locking and stable $n$-cycles. A detailed discussion of this phenomena is in Section \ref{sec:M3bifn}. 

\begin{figure}
\centering
\includegraphics[width=0.5\linewidth]{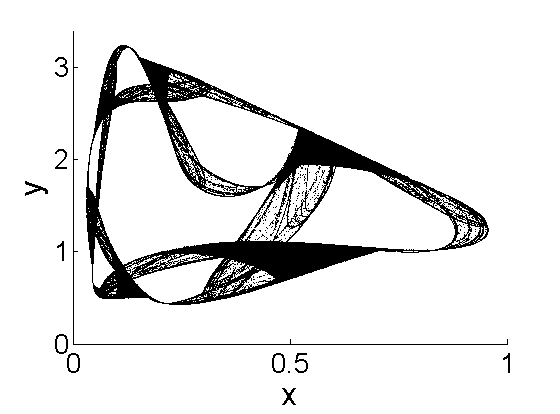}\\
\caption{A strange attractor for Model 4 with $r=2.3$, $b=2.2$. Parameters chosen for aesthetic appeal of the strange attractor. }\label{fig:strange}
\end{figure}

We note that for this model, we also see the development of a chaotic strange attractor. The collapse of the strange attractor in a crisis bifurcation is discussed by Kang et al. \cite{kang2008dynamics}, as well as cases of more complicated bistability between boundary attractors and interior attractors. Hence, we do not discuss details here. One of the strange attractors is shown in Figure \ref{fig:strange}. Due to the use of a stronger nonlinearity in density dependence and stronger parasitism in Model 4, we see the greatest variability in dynamics and bifurcations in the system compared to Models 1, 2, and 3. 

%%%%%%%%%   FIGURE 12 HERE  %%%%%%%%%%%%%%

%%%%%%%%%%%%%%%%%%%%%%%%%%%%%%%%%%%
%%%%%%%      FIGURE 12           %%%%%%%%%%%%%%%%%%%%
%%%%%%%%%%
%%%%%%%%%%%%%%%%%%%%%

%%%%%%%%%%%%%%%%%%%%%%%%%%%%%%%%%%%
\section{Discussion}

We have developed a framework for investigating host--parasitoid systems where density dependence precedes parasitism in the life cycle of the host. Recall that these models have the form given in system (\ref{Eq:GrowthFirst}). Our analysis addresses all combinations of the most frequently used functions for host density-dependence and parasitism. The methods used in this paper can also be extended to models using other functional forms for recruitment and parasitism, including cases where there may not be an explicit expression for the coexistence equilibria. With our analytical approach, we were able to more fully categorize the dynamics of system (\ref{Eq:Model4}), Model 4, which previously had been analyzed using numerical techniques \cite{kang2008dynamics}. Our systematic approach allows for direct comparison of four foundational models, each based on specific biological characteristics of host and parasitoid species. 

Each model resulted in different dynamics. Through systematic comparison of the models, we identified the effects of stronger parasitism (corresponding to higher $\kappa$ or parasitoid aggregation). We then contrasted these effects with the effects of stronger nonlinearity in the density-dependence term. As expected, fractional recruitment and parasitism yield stable dynamics. Stronger parasitism in the model leads to a restricted stability region for the coexistence equilibrium, seen in Figures \ref{fig:Model3StabilityRegion} and \ref{fig:Model4StabilityRegion}. Both Models 3 and 4 include Neimark-Sacker bifurcations where the coexistence equilibrium is replaced with invariant circles. On the other hand, stronger nonlinearity in the density-dependence term produces period-doubling bifurcations and the potential for bistability. In the case of Model 2, the period doubling may be supercritical or subcritical, depending on the value of $b$. The period-doubling bifurcation observed in Model 4 is subcritical and only occurs for sufficiently large values of $r$. 

For models with stronger parasitism resulting from higher parasitoid aggregation (Models 3 and 4), stability of the equilibrium is lost as $b$ increases. Since $b$ is proportional to host carrying-capacity, $K$, an increase in host carrying-capacity can result in loss of stability of the equilibrium for these models, consistent with the paradox of biological enrichment \cite{rosenzweig1971paradox}. For the invariant circles and $n$-cycles that arise after the Neimark--Sacker bifurcation, the host population remains below the carrying capacity throughout the population cycles. On the other hand, the loss of stability through increased $r$ in Model 2 yields drastic swings in host population size above and below carrying capacity, $K$, with relatively short period (2, 4, etc.). In these cases, the introduction of a parasitoid species could increase the host population size about its natural carrying capacity during some years of the population cycles. In agricultural scenarios, these host outbreaks could have devastating consequences.

Future work for the models presented in this paper requires comparison with data from host--parasitoid systems and consideration of what range of parameters are observed biologically. While we have provided a mathematical characterization of these systems, the biological implications need to be experimentally verified. As noted above, the period and amplitude of oscillations differ for the case of invariant circles arising in models with higher parasitism and the case of 2-cycles or 4-cycles arising in models with overcompensatory density-dependent effects. It would be beneficial to compare these models with data to determine if overcompensation does, in fact, lead to shorter-period, higher-amplitude oscillations in host population size in experimental systems. 

In comparing with data, it is important to acknowledge environmental and demographic stochasticity, which will impact the ways that mathematically predicted $n$-cycles and quasi-periodic fluctuations in population size manifest in real populations. It is also important to consider whether the models presented here can be used for prediction in specific management scenarios or whether their use is more suited to development of biological control theory. Barlow \cite{barlow1999models} provides a survey of biological control models for specific real-world systems and emphasizes the value of models in understanding specific case studies, whether or not the models are used for practical management decisions.  

As discussed in Section \ref{Section:Intro}, the sequence of events in the host life-cycle also has important impacts on the population dynamics. The models investigated in this paper assume that density dependence precedes parasitism, which is an appropriate assumption for some species. For example, houseflies (\emph{Musca} spp.) are attacked by pupal parasitoids such as \emph{Spalangia} spp. and \emph{Muscidifurax} spp. after significant density dependence in the early larval stages \cite{may1981density}. However, in other species, density dependence acts on the survivors of parasitism, which leads to the model
\begin{subequations}\label{Eq:Model_Pfirst}
\begin{align}
N_{t+1}&=N_t[1-P_tH(P_t)] G\Big(N_t[1-P_tH(P_t)]\Big),\\
P_{t+1}&=cN_tP_tH(P_t).
\end{align}
\end{subequations}
Note that this model assumes not only that parasitism occurs first in the host life-cycle, but also that the parasitized hosts are functionally dead and unable to compete. For the fractional form for recruitment and the negative binomial form for parasitism, this model can lead to stable equilibria with hosts at a higher level than their carrying capacity in the absence of parasitoids \cite{may1981density, mills1996modelling}.

Systems (\ref{Eq:Model_Pfirst}) and (\ref{Eq:GrowthFirst}) represent the scenarios where parasitism occurs either before or after density-dependent effects on the host, and May et al.\cite{may1981density} compared some specific models in these frameworks. However, more complicated parasitoid phenologies exist in nature. Cobbold et al.\cite{cobbold2009impact} explicitly consider koinobiont parasitoids, which do not kill their host immediately. This means that there is a period of time when parasitized hosts are competing with nonparasitized hosts, which cannot be accounted for with either system (\ref{Eq:Model_Pfirst}) or system (\ref{Eq:GrowthFirst}). Cobbold et al. \cite{cobbold2009impact} found that the delayed mortality of parasitized hosts may have implications for biological control. Differences in the timing of interaction between parasitoids and hosts lead to different predicted population dynamics. Thus, model formulation requires care and awareness of biological assumptions that are inherent to the structure of a model.

Host--parasitoid models have numerous avenues for the inclusion of additional biological complexities such as spatial heterogeneity, Allee effects, and multiple parasitoid species. In building towards these more biologically realistic models, it is important to understand the dynamics of simpler models, such as those analyzed and compared here. Hassell \cite{hassell1978dynamics, hassell2000spatial}  has done excellent work in bridging the gap between simple mechanistic models for host--parasitoid systems and models for more complex and biologically realistic systems. Extending simple mechanistic models to investigate more complicated scenarios can only occur when the simple foundational models are well-understood and presented with explicit acknowledgment of biological assumptions.

Biological differences between models may be critical to communicate well with ecologists and experimentalists. We therefore urge researchers to exercise caution in formulation of models and underlying biological assumptions in order to promote communication and broader understanding of mathematical and theoretical findings.

%%%%%%%%%%%%%%%%%%%%%%%%%%%%%%%%%%%%%%%%%%%%%%%%%%%%
\bibliography{DiscretePredPreybib}
\bibliographystyle{tfs}

%%%%%%%%%%%%%%%%%%%%%%%%%%%%
%%%%%%%%%%%%%%%%%%%%%%%%%%%%%%%%%%%%%%%%%%%%%%%%%%%%%%%%
%% APPENDIX
%%%%%%%%%%%%%%%%%%%%%%%%%%%%%%%%%%%%%%%%%%%%%%%%%%%%%%%%%%%%%%%%%%%%%%%%%%%%%%%%%%%%

\appendix

\section{Partial Derivatives for Jury Conditions}

We begin by evaluating the partial derivatives that appear in the Jury conditions. In doing so, we will use the nullcline equations, $u(x,y)=1$, $v(x,y)=1$. From the definitions of $u(x,y)$ and $v(x,y)$, equations  (\ref{Eq:u}) and (\ref{Eq:v}), we obtain the partial derivatives
\begin{align}
u_x&=g'(x)[1-yh(y)],\\
u_y &=g(x)[1-yh(y)]',\\
v_x &= b[g(x)+xg'(x)]h(y),\\
v_y&=bxg(x)h'(y).
\end{align}

For Models 1 and 2, with fractional parasitism (\ref{Eq:fracxy}), in addition to $u(x,y)=1$ and $v(x,y)=1$, we also use 
\begin{align}
h(y)=\frac{1}{1+y}=1-yh(y). 
\end{align}
Recall that Models 3 and 4 have exponential parasitism, given by equation (\ref{Eq:expxy}). Furthermore, Models 1 and 3 use fractional per-capita-recruitment, with $g(x)$ defined in equation (\ref{Eq:BevHoltxy}), while Models 2 and 4 use exponential per-capita-recruitment, with $g(x)$ defined in equation (\ref{Eq:Rickerxy}). Simplified expressions for the partial derivatives for each model using the corresponding functions for $g(x)$ and $h(y)$ are given in Table \ref{table:pdvs}.

%%%%%%%%%%        TABLE A1 HERE     %%%%%%%%%%%%%%%%%%%%%

%%%%%%%%%%%%%%%%%%%%%%%%%%%%%%%5
%%%%%%%%%%            TABLE A1
%%%%%%%%%%%%%%%%%%%%%%%%%%%%%%%%
\renewcommand{\arraystretch}{2.5}
\begin{table}
\centering
\caption{Partial derivatives used to apply the Jury conditions to the coexistence equilibrium point(s) for each model.}
\begin{tabular}{| c | c | c | c | c |}
\hline
 & Model 1 & Model 2 & Model 3 & Model 4 \\ \hline
\(\displaystyle u_x\) & \(\displaystyle\frac{(1-R_0)}{R_0}g(x)\)                          & \(\displaystyle-r\)                     & \(\displaystyle\frac{(1-R_0)}{R_0}g(x)\)            & \(\displaystyle-r\) \\  \hline
\(\displaystyle u_y\) & \(\displaystyle-h(y)\)                                                          & \(\displaystyle-h(y)\)                 & \(\displaystyle-1\)                                                & \(\displaystyle-1\) \\ \hline 
\(\displaystyle v_x\) &    \(\displaystyle\frac{1}{R_0x}g(x)\)  &  \(\displaystyle\frac{1}{x}-r\)   & \(\displaystyle\frac{1}{R_0x}g(x)\)                    & \(\displaystyle\frac{1}{x}-r\)\\  \hline
\(\displaystyle v_y\) &        \(\displaystyle-h(y)\)                                                   & \(\displaystyle-h(y)\)                 &   \(\displaystyle\frac{1}{yh(y)}[1-yh(y)-h(y)]\)   &    \(\displaystyle\frac{1}{yh(y)}[1-yh(y)-h(y)]\)          \\  \hline
\end{tabular}
\label{table:pdvs}
\end{table} 
\renewcommand{\arraystretch}{1}

%%%%%%%%%%%%%%%%%%%%%%%%%%%%%%%%%%%%%%%%

\section{Model 1 stability calculations}\label{app:M1}
%%%%%%%%%%%%%%%%%%%%%%%%%%%%%%%%%%%%
\subsection{Requirements for existence of coexistence equilibrium in the first quadrant}\label{app:Model1FirstQuadrant}

We now determine the conditions that ensure an equilibrium in the interior of the first quadrant. For this model, we can explicitly solve system (\ref{Sys:uv1}) for the coexistence equilibrium,
\begin{align}\label{Model1Equil}
(x^*,y^*)&=\left(\frac{1}{b},g\left(\frac{1}{b}\right)-1\right).
\end{align}
The $x$ coordinate is positive for all positive $b$. The $y$ coordinate is positive when
\begin{align}
g\left(\frac{1}{b}\right)-1&=\frac{R_0}{1+\left(R_0-1\right)\left(\frac{1}{b}\right)}-1>0,
\end{align}
which simplifies to $b>1$ since we assume $R_0>1$. Thus, the coexistence equilibrium exists and is in the first quadrant when $R_0>1, b>1$. 
%
%%%%%%%%%%%%%%%%%%%%%%%%%%%%%%
\subsection{First Jury condition}\label{app:Model1Jury1}
For $b>1$, $R_0>1$, the $x$ and $y$ coordinates of the coexistence equilibrium are positive. We thus use partial derivatives from Table \ref{table:pdvs} to write inequality (\ref{Jury1part}), as
\begin{align}
u_xv_y-u_yv_x=\left(\frac{R_0-1}{R_0}\right)g(x)h(y)+h(y)\left[\frac{1}{R_0x}g(x)\right]&>0,
\end{align}
which simplifies to
\begin{align}
\frac{h(y)}{x}&>0.
\end{align}
Since $h(y)$ is positive, the first Jury condition is satisfied whenever the coexistence equilibrium is in the first quadrant. 

When $b=1$, the $y$-coefficient from equation (\ref{Model1Equil}) is $y^*=0$, and the first Jury condition, inequality (\ref{Jury1})  is violated. When $R_0=1$, equation (\ref{Model1Equil}) again gives us $y^*=0$, regardless of the value of $b$, such that the $x$-axis  is a line of equilibrium points. For $R_0=1$, the first Jury condition, inequality (\ref{Jury1}) is again violated.

\subsection{Second Jury condition}\label{app:Model1Jury2}
Recall that the second Jury condition, inequality (\ref{Eq:Jury2}), is
\begin{align}
1+\tau+\Delta& >0.
\end{align}
For this model, we will not show this directly. Instead, note that if $\tau>0$ and $1-\tau+\Delta>0$, which is the first Jury condition, then $1+\tau+\Delta>1-\tau+\Delta>0$. This means $\tau>0$ and the satisfaction of the first Jury condition are sufficient criteria for the second Jury condition.

The first Jury condition is satisfied for $b>1, R_0>1$. We will show that in this case, the second Jury condition will also be satisfied. We proceed by showing that $\tau>0$ at the equilibrium. As seen in equation (\ref{Eq:tau}), $\tau = 2 +xu_x+yv_y$. We use the expressions for $u_x$ and $v_y$ from  Table \ref{table:pdvs} and the definitions of $g(x)$ and $h(y)$ from equations (\ref{Eq:BevHoltxy}) and (\ref{Eq:fracxy}) to express the trace, 
\begin{equation}\label{Eq:tauJury2}
\tau=2-\frac{R_0-1}{R_0}xg(x)-yh(y)=2-\left[\frac{(R_0-1)x}{1+(R_0-1)x}+\frac{y}{1+y}\right].
\end{equation}
We thus seek to show that
\begin{equation}\label{Eq:Jury2ShowPos}
2>\frac{(R_0-1)x}{1+(R_0-1)x}+\frac{y}{1+y}
\end{equation}
for $x,y>0$, $R_0>1$. 

Both of the terms on the right-hand side of inequality (\ref{Eq:Jury2ShowPos}) are of the form $z(1+z)^{-1}$, where $z$ is positive. Each term individually is less than one because $z<1+z$, which indicates that $z(1+z)^{-1}<1$ for positive $z$.
Therefore,
\begin{align}
\tau=2-\left[\frac{(R_0-1)x}{1+(R_0-1)x}+\frac{y}{1+y}\right]>0.
\end{align}
It follows that the first Jury condition is a sufficient condition for the second Jury condition for Model 1.

For either $b=1$ or $R_0=1$, we can directly calculate the terms in the second Jury condition, inequality (\ref{Eq:Jury2}). Direct calculation verifies that the second Jury condition is satisfied. 

\subsection{Third Jury condition}\label{app:Model1Jury3}
Recall that the third Jury condition is $\Delta<1$.
The determinant is given in terms of the partial derivatives in equation (\ref{DetExp}). Using the expressions for $u_x,u_y, v_x$, and $v_y$ from Table \ref{table:pdvs}, the third Jury condition
simplifies to
\begin{equation}
1+\frac{(1-R_0)}{R_0}xg(x)<1.
\end{equation}
When we substitute equation (\ref{Eq:BevHoltxy}) for $g(x)$, the condition can be expressed as
\begin{equation}
1-\frac{(R_0-1)x}{1+(R_0-1)x}<1,\\
\end{equation}
which simplifies to
\begin{equation}
1+(R_0-1)x>1.
\end{equation}
This is true for $R_0>1$ for the equilibrium in the interior of the first quadrant. Thus, the third Jury condition is satisfied for $R_0>1$. For $R_0=1$, the third Jury condition is violated.
%%%%%%%%%%%%%%%%%%%%%%%%%%%%%%%%%%%%%%%%%%%%
%%%%%%%%%%%%%%%%%%%%%%%%%%%%%%%%%%%%%
\section{Model 2 stability calculations}\label{app:M2}

%%%%%%%%%%%%%%%%%%%%%%%%%%%%%%%%%%%%%%%%%%%%%%%%%%%%%%%%%%%%%%%%%%%%%%%%%%%%%%%%%%%%%%%%%%%%%%%%%%
\subsection{Requirements for existence of coexistence equilibrium in the first quadrant}\label{app:Model2FirstQuadrant}
For this model, we can again explicitly solve system (\ref{Sys:uv1}) for the coexistence equilibrium for Model 2,
\begin{equation}\label{M2Equil}
(x^*,y^*)=\left(\frac{1}{b}, e^{r\left(1-1/b\right)}-1\right).
\end{equation}
When $b=1$, this equilibrium point is on the $x$-axis at $(x^*,y^*)=(1,0)$, which is the exclusion equilibrium. For the coexistence equilibrium to be in the interior of the first quadrant, it is necessary that 
\begin{equation}
y^*=e^{r\left(1-1/b\right)}-1>0,
\end{equation}
For $r>0$, this requires $b>1$. Note that we will not consider the case $r<0$, $b<1$ since we are interested in cases where the host species persists in the absence of the parasitoid.
%%%%%%%%%%%%%%%%%%%%%%%%%%%%%%%%%%%%%%%%%%%%%%%%%%%
\subsection{First Jury condition: slopes of zero-growth isoclines}\label{app:Model2Jury1}
Using partial derivatives from Table \ref{table:pdvs}, the first Jury condition, inequality (\ref{Jury1}), is
\begin{align}
xy\left[rh(y)+h(y)\left(\frac{1}{x}-r\right)\right]=yh(y)&>0,
\end{align}
Because $yh(y)$ is positive for the coexistence equilibrium, this inequality holds for the equilibrium in the interior of the first quadrant. When $b=1$, $yh(y)=0$, and the first Jury condition is violated. For $r=0$, the $x$-axis  is a line of equilibrium points, and the first Jury condition is again violated. 
%%%%%%%%%%%%%%%%%%%%%%%%%%%%%%%%%%%%%%%%%%%%%%%%%%%%%
\subsection{Second Jury Condition} \label{app:Model2Jury2}

Again using partial derivatives from Table \ref{table:pdvs}, the second Jury condition, inequality (\ref{Jury2}) simplifies to
\begin{equation}
4-2xr-yh(y)=4-2xr-\frac{y}{1+y}>0.
\end{equation}
The coordinates of the coexistence equilibrium point are given by equation (\ref{M2Equil}). Using these values, the stability condition is
\begin{equation}\label{app:M2J2}
3-\frac{2}{b}r+e^{\left(\frac{r}{b}-r\right)}>0.
\end{equation}

We now consider the transcendental equation, 
\begin{equation}
3-\frac{2}{b}r+e^{\left(\frac{r}{b}-r\right)}=0,
\end{equation}
and introduce the parameter $u=r/b$ so that
\begin{equation}
3-2u+e^{u-r}=0.
\end{equation}
We solve for $r$ as a function of $u$,
\begin{equation}\label{M2Jury2_r}
r=u-\ln\left(2u-3\right),
\end{equation}
and can then also write $b$ as a function of $u$,
\begin{equation}\label{M2Jury2_b}
b=\frac{r}{u}=1-\frac{1}{u}\ln\left(2u-3\right).
\end{equation}
For $u>3/2$, equations (\ref{M2Jury2_r}) and (\ref{M2Jury2_b}) express the boundary of the region in parameter space where the coexistence equilibrium satisfies the second Jury condition.

For $b=1$, inequality (\ref{app:M2J2}) requires $r<2$. The point $(r,b)=(2,1)$ is where the Jury 2 curve intersects the $b=1$ line, seen in Figure \ref{fig:Model2StabilityRegion}. 
%%%%%%%%%%%%%%%%%%%%%%%%%%%%%%%%%%%%%%%%%%%%%%%%%%%
\subsection{Third Jury Condition}\label{app:Model2Jury3}
The expression for the determinant from equation (\ref{DetExp}) for this model simplies significantly to  
\begin{equation}
\Delta = 1-rx,
\end{equation}
using the partial derivatives in Table \ref{table:pdvs}. Since $x=1/b$ at the equilibrium, the third Jury condition is
\begin{equation}
1-\frac{r}{b}<1.
\end{equation}
Since $b>0$ and we assumed $r\geq 0$, this inequality is satisfied for $r>0$. When $r=0$, the third Jury condition is violated.
%%%%%%%%%%%%%%%%%%%%%%%%%%%%%%%%%%%%%%%%%%%%%%%%%%%%%%%%%%%%%%%%%%%
%%% MODEL 3 APPENDIX %%%
%%%%%%%%%%%%%%%
\section{Model 3 stability calculations}\label{app:M3}
%%%%%%%%%%%%%%%%%%%%%%%%%%%%%%
\subsection{Requirements for existence of coexistence equilibrium in first quadrant}\label{app:Model3FirstQuadrant}
As was true in Section \ref{app:Model1FirstQuadrant}, we seek to determine the conditions that ensure that an equilibrium exists in the interior of the first quadrant, this time for Model 3, system (\ref{Eq:Model3}). The coexistence equilibrium cannot be solved for explicitly in this case, so we instead consider the nullclines.

Equation (\ref{Eq:uis1}) is the host nullcline with intercepts $(0,\ln R_0)$ and $(1,0)$. To obtain the slope of this nullcline in the $x$-$y$ plane, we first differentiate $u(x,y)=1$ and get
\begin{align}
u_x+u_y\frac{dy}{dx}&=0.
\end{align}
The slope of the host nullcline is
\begin{equation}
\frac{dy}{dx}=-\frac{u_x}{u_y}
=\frac{1-R_0}{R_0}g(x),
\end{equation}
using the expressions for $u_x$ and $u_y$ from Table \ref{table:pdvs}. Since $g(x)>0$ and we assume $R_0>1$, the host nullcline is monotone decreasing in the first quadrant from $(0,\ln R_0)$ to $(1,0)$.

We now consider the parasitoid nullcline, equation (\ref{Eq:vis1}). To find the slope in the $x$-$y$ plane, we differentiate $v(x,y)=1$ with respect to $x$ to get
\begin{align}
v_x+v_y\frac{dy}{dx}&=0.
\end{align}
 The slope for the parasitoid nullcline is thus
\begin{equation}
\frac{dy}{dx}=\frac{-v_x}{v_y}=-\frac{g(x)}{R_0xv_y}.
\end{equation}
Since $g(x)>0$, the sign of $v_y$ will determine the sign of the slope of the parasitoid nullcline. Negative $v_y$ will indicate that the slope of the nullcline is positive.

We substitute $h(y)$ from equation (\ref{Eq:expxy}) into $v_y$ for Model 3, such that
\begin{equation}
v_y=\frac{1}{1-e^{-y}}\left[e^{-y}-\frac{1}{y}\left(1-e^{-y}\right)\right]=\frac{1}{y\left(1-e^{-y}\right)}\left(ye^{-y}-1+e^{-y}\right).
\end{equation}
The denominator is positive for $y>0$, so we consider the numerator. For $y>0$,
\begin{subequations}
\begin{align}
1+y&<e^y,\\
(1+y)e^{-y}&<1,\\
e^{-y}+ye^{-y}-1&<0.
\end{align}
\end{subequations}
Thus, we conclude that $v_y<0$ for $y>0$. This means that the slope of the parasitoid nullcline is positive in the first quadrant. If there is an intersection of the host and parasitoid nullclines in the first-quadrant, it is unique.

To determine existence of the equilibrium, we examine the the $x$- and $y$-intercepts of the parasitoid nullcline,
\begin{equation}
1=\frac{bxR_0h(y)}{1+(R_0-1)x}.
\end{equation}
After solving for $x$, we obtain
\begin{equation}\label{M3Par_xint}
x=\frac{1}{bR_0h(y)+(1-R_0)},
\end{equation}
where
\begin{equation}
h(y)=\frac{1}{y}\left(1-e^{-y}\right).
\end{equation}
To examine equation (\ref{M3Par_xint}), we consider the limiting behavior of $h(y)$ as $y\to-\infty$,
\begin{equation}
\lim_{y\to-\infty}h(y)=\lim_{z\to\infty}h(-z)=\lim_{z\to\infty} \frac{-1}{z}\left(1-e^z\right)=\infty.\\
\end{equation}
Thus, if we consider the limit as $y\to-\infty$ in equation (\ref{M3Par_xint}), $x\to 0^+$. This nullcline does not have a $y$-intercept because as $x\to 0^+$, $y\to-\infty$.

Since we know that in the first quadrant, the parasitoid nullcline is monotone increasing and the host nullcline has $x$-intercept at $x=1$, we need to find the conditions for which the parasitoid nullcline's $x$-intercept lies between 0 and 1. For these conditions, there exists exactly one intersection of the parasitoid and host nullclines in the interior of the first quadrant. The $x$-intercept of the parasitoid nullcline is the solution to
\begin{equation}
1=bxg(x)h(0)=\frac{bxR_0}{1+(R_0-1)x},
\end{equation}
which is
\begin{align}
x_{\text{int}}=\frac{1}{R_0(b-1)+1}.
\end{align}

We seek the conditions for which 
\begin{align}
0&<\frac{1}{R_0(b-1)+1}<1.
\end{align}
This translates into the following two criteria,
\begin{align}
R_0(b-1)+1&>0,
\end{align}
and
\begin{align}
R_0(b-1)+1 &>1,
\end{align}
which can be consolidated as
\begin{equation}
\begin{split}
R_0(b-1)+1&>1,\\
R_0(b-1)&>0.
\end{split}
\end{equation}
This is true when $b>1$. So the $x$-intercept of the parasitoid nullcline occurs between 0 and 1 if and only if $b>1$. 

We conclude that there is exactly one equilibrium point in the interior of the first quadrant if and only if $b>1$. When $b>1$, we can then determine if the coexistence equilibrium is stable. For $b=1$, the equilibrium point is on the boundary of the first quadrant, at $(1,0)$. For $R_0>1$, $b<1$, there are no equilibria points in the interior of the first quadrant.
%
%%%%%%%%%%%%%%%%%%%%%%%%%%%
\subsection{First Jury condition}\label{app:Model3Jury1}

For $b>1$, $R_0>1$, the $x$ and $y$ coordinates of the coexistence equilibrium are positive. We thus use partial derivatives from Table \ref{table:pdvs} to write the first Jury condition, inequality (\ref{Jury1part}), as
\begin{align}
xy(u_xv_y-u_yv_x)=xy\left[\frac{1-R_0}{R_0}g(x) v_y +\frac{1}{R_0x} g(x)\right] &>0,
\end{align}
which simplifies to
\begin{align}\label{M3J1}
\frac{1}{R_0x} g(x) &>\frac{R_0-1}{R_0}g(x) v_y.
\end{align}
The left-hand side of inequality (\ref{M3J1}) is positive, while the right-hand side is negative for $R_0>1$, since $v_y$ was shown to be negative in Section \ref{app:Model3FirstQuadrant}. Thus, this inequality holds for the positive coexistence equilibrium.

When $b=1$, the coexistence equilibrium has collided with the exclusion equilibrium at $(1,0)$. Since the $y$ coordinate is $0$, the first Jury condition, inequality (\ref{Jury1}),  is violated. When $R_0=1$, any point on the $x$-axis is a solution to system (\ref{Eq:Model3}). Since these equilibria points have $y=0$, the first Jury condition is violated for this line of equilibrium points.

%%%%%%%%%%%%%%%%%%%%%%%%%%%%%%%%%%%%%
\subsection{Second Jury condition}\label{app:Model3Jury2}
We will use the technique from Section \ref{app:Model1Jury2} to show that the first Jury condition is a sufficient condition for the second Jury condition. To do this, we must show that $\tau>0$ at the interior equilibrium. As seen in equation (\ref{Eq:tau}), $\tau = 2 +xu_x+yv_y$. We use the expressions for $u_x$ and $v_y$ from Table \ref{table:pdvs} and the definitions of $g(x)$ and $h(y)$ from equations (\ref{Eq:BevHoltxy}) and (\ref{Eq:expxy}) to get
\begin{equation}
\tau=1-\left[\frac{(R_0-1)x}{1+(R_0-1)x}\right]+\frac{ye^{-y}}{1-e^{-y}},
\end{equation}
after simplification. 

For positive $y$, the last term is positive. Similarly to Section \ref{app:Model1Jury2}, the middle term is of the form $z(1+z)^{-1}$, where $z$ is positive. For $R_0>1$, this term individually is less than one because $z<1+z$, which indicates that $z(1+z)^{-1}<1$ for positive $z$. Since the coordinates of the coexistence equilibrium are positive for $b>1$, $R_0>1$, we conclude that $\tau>0$ for $b>1$ and $R_0>1$. It follows that the first Jury condition is a sufficient condition for the second Jury condition.

For either $b=1$ or $R_0=1$, we can directly calculate the terms in the second Jury condition, inequality (\ref{Eq:Jury2}). Direct calculation verifies that the second Jury condition is satisfied in these cases. 
%%%%%%%%%%%%%%%%%%%%%%%%%%%%%%%%%%%%%%%
\subsection{Third Jury condition}\label{app:Model3Jury3}
The third Jury condition is $\Delta <1$. The determinant is given in terms of the partial derivatives in equation (\ref{DetExp}). Using the expressions for $u_x,u_y, v_x$, and $v_y$ from Table \ref{table:pdvs} and much algebraic simplification, the third Jury condition is
\begin{equation} \label{M3J3}
\Delta=\frac{g(x)}{R_0h(y)}<1.
\end{equation}

We want to write the condition solely in terms of the parameters, $R_0$ and $b$, and find the curve in the $R_0$-$b$ plane where stability changes. Because of the transcendental nature of the inequality, we will express this curve parametrically with $R_0$ and $b$ as functions of $y$. To do so, we consider
\begin{equation}
1=\frac{g(x)}{R_0h(y)}=\frac{1}{h(y)\left[1+(R_0-1)x\right]},\label{M3J3_eq}
\end{equation}
and solve for $x$,
\begin{equation}
x=\left(\frac{1}{R_0-1}\right)\left[\frac{1-h(y)}{h(y)}\right].\label{M3J3_x}
\end{equation}

We now incorporate the host nullcline, equation (\ref{Eq:uis1}), which is valid at the equilibrium point. Using $g(x)=R_0h(y)$ from equation (\ref{M3J3_eq}), we get
\begin{equation}
1=g(x)[1-yh(y)]=R_0h(y)[1-yh(y)].
\end{equation}
Solving for $R_0$ as a function of $y$ yields
\begin{equation}\label{Eq:M3J3_R0}
R_0=\frac{1}{h(y)[1-yh(y)]}.
\end{equation}

Next, we need an expression for $b$ as a function of $y$. To do this, we incorporate the parasitoid nullcline, equation (\ref{Eq:vis1}), which is valid at the equilibrium point. Starting with equation (\ref{Eq:vis1}), we replace $x$ with the expression from (\ref{M3J3_x}) and also substitute $R_0h(y)$ for $g(x)$, using the determinant condition (\ref{M3J3_eq}). This gives us the equation,
\begin{equation}
1=bxg(x)h(y)=
\frac{b}{R_0-1}\left[\frac{1-h(y)}{h(y)}\right]R_0h(y)h(y).
\end{equation}
We solve for $b$ to get
\begin{equation}\label{Eq:M3J3_b_R0}
b=\left(\frac{R_0-1}{R_0}\right)\frac{1}{h(y)\left[1-h(y)\right]}.
\end{equation}

We then eliminate the dependence on $R_0$ from the equation for $b$. This gives us $b$ as a function of $y$,
\begin{align}
b&=\frac{1-h(y)[1-yh(y)]}{h(y)[1-h(y)]},
\end{align}
which does not simplify in a meaningful way. This equation combined with equation (\ref{Eq:M3J3_R0}) expresses the boundary of the region in parameter space where the coexistence equilibrium satisfies the third Jury condition.

When $R_0=1$, the $x$-axis is a line of equilibrium points, as stated in Section \ref{app:Model3Jury1}. Under these conditions, the expression for the determinant simplifies to $\Delta =1$, and so the third Jury condition is also violated for $R_0=1$.
%%%%%%%%%%%%%%%%%%%%%%%%%%%%%%%%%%%%%%%%%%%%%%%%%%%%%%%%%%%%%%%%%%%%%%%%%%%%%%%%%%%%%%%%%%%%%%%%%%%%%%%%%%%%%%%%%%%%%%%%%%%%%

\section{Model 4 stability calculations}\label{app:M4}
From Section \ref{Section:M4}, recall that there may be one or two coexistence equilibria, depending on the parameter values. In the case of two coexistence equilibria, only the point with the larger $y$ value may be stable, as discussed in Section \ref{Section:M41}. The analysis here pertains to the stability of the single unique coexistence equilibrium or the coexistence equilibrium point with the larger $y$ value. 

As Kang et al. \cite{kang2008dynamics} proved, for $b>1$, system (\ref{Eq:Model4}) has a unique positive equilibrium. For $b=1$, $0<r<2$, the point $(1,0)$ is an equilibrium point, and there is no coexistence equilibrium in the interior of the first quadrant. For $r>2$ and $b$ just less than $1$, the system has both a stable coexistence equilibrium point and an unstable coexistence equilibrium point, as seen in Figure \ref{fig:AlmostCollision}. For $r>2$ and $b=1$, the unstable coexistence point collides with the exclusion equilibrium, $(1,0)$. This is all consistent with the analysis in Kang et al. \cite{kang2008dynamics}.

%%%%%%%%%%%%%%%%%%%%%%%%%%%%%%%%%%%%%%%%%%%%%%%%%%%
\subsection{First Jury condition}\label{app:Model4Jury1}
For $b>1$, $r>0$, the $x$ and $y$ coordinates of the coexistence equilibrium are positive. Using partial derivatives from Table \ref{table:pdvs}, the first Jury condition, inequality (\ref{Jury1}), is
\begin{align}
xy\left\{-r\left[\frac{1}{yh(y)}\right]\left[1-yh(y)-h(y)\right]+\frac{1}{x}-r\right\} &>0,
\end{align}
which simplifies to
\begin{align}\label{M4J1}
r\left(\frac{e^y-1-ye^y}{ye^y-y}\right)+\frac{1}{x}&>0,
\end{align}
using equation (\ref{Eq:expxy}) for $h(y)$. We want to write the condition solely in terms of the parameters, $r$ and $b$, and find the curve in the $r$-$b$ plane where stability changes. Because of the transcendental nature of the inequality, we will write this curve parametrically with $r$ and $b$ as functions of $y$. To do so, we first consider
\begin{align}\label{M4J1eq}
r\left(\frac{e^y-1-ye^y}{ye^y-y}\right)+\frac{1}{x}&=0.
\end{align}

We now incorporate the host and parasitoid nullclines, equations (\ref{Eq:uis1}) and (\ref{Eq:vis1}). For Model 4, these equations simplify to 
\begin{align}\label{M4uis1}
r-rx-y&=0,
\end{align}
and
\begin{align}\label{M4vis1}
bxe^{r-rx}\ \frac{1}{y}\left(1-e^{-y}\right)&=1.
\end{align}
Using equation (\ref{M4uis1}), we now eliminate $x$ from equation (\ref{M4J1eq}) and write $r$ as a function of $y$, which simplifies to
\begin{align}\label{M4J1r}
r&=\frac{y^2e^y}{1+ye^y-e^y}.
\end{align}

We now return to equation (\ref{M4vis1}) and again eliminate $x$. We can then write $b$ as a function of $y$, using equation (\ref{M4J1r}) to eliminate $r$. After algebraic simplification, we obtain
\begin{align}\label{M4J1b}
b&=\frac{y^2e^y}{(e^y-1)^2}.
\end{align}
Equations (\ref{M4J1r}) and (\ref{M4J1b}) give the boundary of the region in parameter space where the coexistence equilibrium satisfies the first Jury condition. Since this curve is just barely below the curve for the second Jury condition found in Section \ref{app:Model4Jury2}, this curve does not contribute to the stability region shown in Figure \ref{fig:Model4StabilityRegion}.

Now consider $b=1$ with $0<r<2$. As discussed previously, there is no equilibrium in the interior of the first quadrant for this case because the coexistence equilibrium has collided with the exclusion equilibrium point at $(1,0)$. Since $y=0$, the first Jury condition is violated for $b=1$, $0<r<2$. Note also that when we take the limit as $y\to 0$ for equations (\ref{M4J1r}) and (\ref{M4J1b}), we obtain $(r,b) = (2,1)$. This means that the first Jury condition curve described by equations (\ref{M4J1r}) and (\ref{M4J1b}) connects to the first Jury condition curve given by $b=1$, $0<r<2$. Finally, when $r=0$, system (\ref{Eq:Model4}) has a line of equilibria on the $x$-axis. Since each of these equilibrium points has $y=0$, the first Jury curve is violated for $r=0$.

%%%%%%%%%%%%%%%%%%%%%%%%%%%%%%%%%%%%%%%%%%%%%%%%%%%%%
\subsection{Second Jury Condition} \label{app:Model4Jury2}
We use partial derivatives from Table \ref{table:pdvs} and equation (\ref{Eq:expxy}) to write the second Jury condition, inequality (\ref{Jury2}), as
\begin{align}\label{M4J2}
4+2x(-r)+2y\left[\frac{y-e^y+1}{y(e^y-1)}\right]+xy\left[\frac{-r (y-e^y+1)}{y(e^y-1)}\right]-xy\left[-1\left(\frac{1}{x}-r\right)\right]&>0.
\end{align}
Our goal is now to determine a curve in parameter space where stability changes. We re-write inequality (\ref{M4J2}) as a equality and eliminate $x$ using the expression from the host nullcline given in equation (\ref{M4uis1}). When we simplify and solve for $r$, we obtain
\begin{align}\label{M4J2r}
r&=\frac{2e^y-2+2ye^y+y^2e^y}{e^y-1+ye^y}.
\end{align}

To get $b$ as a function of $y$, we first use equation (\ref{M4uis1}) to eliminate $x$ from the parasitoid nullcline, equation (\ref{M4vis1}). Then, we use equation (\ref{M4J2r}) to eliminate $r$. We solve for $b$ and obtain
\begin{align}\label{M4J2b}
b&=\frac{2y(e^y-1)+y^2e^y(2+y)}{(2+y)e^{2y}-4e^y+2-y}.
\end{align}

Equations (\ref{M4J2r}) and (\ref{M4J2b}) give the boundary of the region in parameter space where the coexistence equilibrium satisfies the second Jury condition. Note that when we take the limit as $y\to 0$ for equations (\ref{M4J2r}) and (\ref{M4J2b}), we obtain $(r,b) = (2,1)$. For $b=1$, inequality (\ref{M4J2}) requires $r<2$. The point $(r,b)=(2,1)$ is where the Jury 2 curve intersects the $b=1$ line. Thus, at the point $(r,b)=(2,1)$, the second Jury curve given parametrically by equations (\ref{M4J2r}) and (\ref{M4J2b}) intersects the first Jury curve, which is described in section \ref{app:Model4Jury1}.

%%%%%%%%%%%%%%%%%%%%%%%%%%%%%%%%%%%%%%%%%%%%%%%%%%%
\subsection{Third Jury Condition}\label{app:Model4Jury3}
We use partial derivatives from Table \ref{table:pdvs} to write the third Jury condition, inequality (\ref{Jury3}), as
\begin{align}
1-rx+ \frac{y}{yh(y)}[1-yh(y)-h(y)]-\frac{rxy}{yh(y)}[1-yh(y)-h(y)]+xy\left(\frac{1}{x}-r\right)&<1.
\end{align}
This simplifies to 
\begin{align}
\frac{ye^y}{e^y-1}(1-xr)&<1,
\end{align}
where we use equation (\ref{Eq:expxy}) for $h(y)$. 

To find the curve where stability of the equilibrium changes, we again consider an equation instead of the inequality. After eliminating $x$ using equation (\ref{M4uis1}) from the host isocline, we solve for $r$ as a function of $y$,
\begin{align}\label{M4J3r}
r&=\frac{ye^y+y^2e^y-e^y+1}{ye^y}.
\end{align} 

As was done in Section \ref{app:Model4Jury2}, we use the parasitoid nullcline, equation (\ref{M4vis1}), with equation (\ref{M4J3r})
to write $b$ as a function of $y$,
\begin{align}\label{M4J3b}
b&=\frac{y^3e^y+y^2e^y-ye^y+y}{(e^y-1)(ye^y-e^y+1)}.
\end{align}
Equations (\ref{M4J3r}) and (\ref{M4J3b}) give the boundary of the region in parameter space where the coexistence equilibrium satisfies the third Jury condition.

Additionally, when $r=0$, the $x$-axis is a line of equilibrium points, as stated in Section \ref{app:Model4Jury1}. Under these conditions, the expression for the determinant simplifies to $\Delta =1$, and so the third Jury condition is also violated for $r=0$ ($R_0=1$). 

%%%%%%%%%%%%%%%%%%%%%%%%%%%%%%%%%%%%%%%%%%%%%%%%%%%%%%

\end{document}